\PassOptionsToPackage{usenames,dvipsnames,table}{xcolor}
\documentclass[USenglish,10pt,a4paper,arxiv]{article}

	\PassOptionsToPackage{noend}{algpseudocode}
	\PassOptionsToPackage{capitalize,nameinlink}{cleveref}
	\PassOptionsToPackage{breakspaces}{clevethm}
	\PassOptionsToPackage{dvipsnames}{xcolor}

	\usepackage[cal=boondoxo, scr=euler]{mathalfa}
	\usepackage[bottom=1.2in,left=1.1in,right=1.1in]{geometry}
	\usepackage[
		alglinenumber,
		eqreset=section,
		thmreset=section,
	]{myPreamble}

	\let\OLDand\and
	\def\and{\texorpdfstring{\OLDand}{, }}%
	\newenvironment{keywords}{\par\noindent{\bf Keywords. }}{}
	\newenvironment{AMS}{\par\noindent{\bf AMS subject classifications. }}{}
	\newcommand{\Sep}{ \(\cdot\)\ }

	\newenvironment{assumption}{\begin{ass}}{\end{ass}}
	\newenvironment{corollary}{\begin{cor}}{\end{cor}}
	\newenvironment{definition}{\begin{defin}}{\end{defin}}
	
	\newenvironment{lemma}{\begin{lem}}{\end{lem}}
	
	\newenvironment{remark}{\begin{rem}}{\end{rem}}
	\newenvironment{theorem}{\begin{thm}}{\end{thm}}

	\allowdisplaybreaks

	\makeatletter
		\renewcommand{\clevethm@proofsectiontitle}{Omitted proofs of }
	\makeatother
	\Crefname{assumption}{Assumption}{Assumptions}
	\Crefname{example}{Example}{Examples}
	\Crefname{fact}{Fact}{Facts}
	\Crefname{remark}{Remark}{Remarks}
	\crefname{ALG@line}{step}{steps}
	\crefname{enumeratpropi}{property}{properties}
	\crefname{enumeratpropii}{property}{properties}
	\counterwithout{algorithm}{section}
	\renewcommand{\thealgorithm}{\arabic{algorithm}}
	\crefformat{equation}{(#2#1#3)}

	\DeclareMathAlphabet{\mathcalligra}{T1}{calligra}{m}{n}
	\newcommand{\C}{\mathcal C}							
	\newcommand{\LL}{\mathcal L}						
	\newcommand{\LM}{{\text{\sc lm}}}					
	\newcommand{\M}{\mathcal M}							
	\newcommand{\s}{\tilde s}
	\newcommand{\ki}[1][i]{t_\nu(#1)}
	\newcommand{\U}{\mathcal U}

	\newcommand{\Op}{\mathcal D}						
	\newcommand{\loss}{\mathscr L}		 				

\makeatletter
	\def\operator@font{\rm}
	
	\DeclareMathOperator{\BregmanD}{D}
	\let\mod\relax
	\DeclareMathOperator{\mod}{mod}
	\DeclareMathOperator{\poisson}{Poisson}%

	\newcommand{\FBE}{\Phi^{\h*}}						
	\newcommand{\FBEG}{\psi^{\h}}						
	\newcommand{\T}{\operatorname T_{\Gamma}^{F,G}}		

	\newcommand{\BregmanKernel}{H}						
	\newcommand{\D}{\@ifstar\@@D\@D}					
	\newcommand{\h}{\@ifstar{\@@h}{\@h}}				
	\let\P\relax
	\newcommand{\E}{\@ifstar\@@E\@E}					
	\newcommand{\P}{\@ifstar\@@P\@P}					

	\newcommand{\@h}{\@ifnextchar_{h}{\BregmanKernel}}
	\newcommand{\@@h}{\@ifnextchar_{\hat h}{\hat\BregmanKernel}}

	\newcommand{\@D}{\@ifnextchar_{\@D@sub}{\BregmanD_{\@h}}}
	\newcommand{\@@D}{\@ifnextchar_{\@@D@sub}{\@@Dnosub}}
	\def\@D@sub_#1{\BregmanD_{\@h_{#1}}}
	\def\@@D@sub_#1{\BregmanD_{\@@h_{#1}}}

	\newcommand{\@@Dnosub}{\@ifstar{\@@D@@}{\BregmanD_{\@@h}}}
	\newcommand{\@@D@@}{\@ifnextchar_{\@@D@@sub@}{\BregmanD_{\@@h@}}}
	\newcommand{\@@h@}{\@ifnextchar_{\conj{\hat h}}{\conj{\hat\BregmanKernel}}}
	\def\@@D@@sub@_#1{\BregmanD_{\@@h@_{#1}}}
	
	\newcommand{\@P}[2][k]{\mathcal P_{#1}\ifstrempty{#2}{}{\left[#2\right]}}
	\newcommand{\@@P}[2][k]{\mathcal P_{#1}\ifstrempty{#2}{}{[#2]}}

	\newcommandx{\@E}[3][1=k,3={}]{
		\mathbb E_{#1}\ifstrempty{#2}{}{
			\left[
				#2\ifstrempty{#3}{}{\vphantom{#3}\right|\left.#3}
			\right]
		}
	}
	\newcommandx{\@@E}[3][1=k,3={}]{
		\mathbb E_{#1}\ifstrempty{#2}{}{
			[#2\ifstrempty{#3}{}{|#3}]
		}
	}

	\newcommand{\bl}{\@ifstar\@@bl\@bl}
	\newcommand{\@bl}[1]{#1}
	\newcommand{\@@bl}[1]{#1}
\makeatother

	\newcommand{\ralign}[1]{\fillwidthof[l]{\ralignToWidthOf}{#1}}%
	\newcommand{\ralignToWidthOf}{}

	\pgfplotsset{compat=1.16}

\showgrayfalse

	\newcommand{\TheTitle}{Bregman Finito/MISO for nonconvex regularized finite sum minimization without Lipschitz gradient continuity}
	\newcommand{\TheShortTitle}{Bregman Finito/MISO for nonconvex finite sum minimization}
	
		\newcommand{\TheFunding}{%
		This work was supported by the Research Foundation Flanders (FWO) PhD grant 1196820N and research projects G0A0920N, G086518N, and G086318N;
		Research Council KU Leuven C1 project No. C14/18/068;
		Fonds de la Recherche Scientifique -- FNRS and the Fonds Wetenschappelijk Onderzoek -- Vlaanderen under EOS project no 30468160 (SeLMA);
		JSPS KAKENHI grant number JP21K17710.%
	}
	\newcommand{\TheKeywords}{%
		Nonsmooth nonconvex optimization\Sep
		incremental aggregated algorithms\Sep
		Bregman Moreau envelope\Sep
		KL inequality%
	}
	\newcommand{\TheAMSsubj}{%
		90C06\Sep 
		90C25\Sep 
		90C26\Sep 
		49J52\Sep 
		49J53.
	}
	\newcommand{\TheAddressUA}{%
		Department of Mathematics and Computer Science, University of Antwerp,
		Middelheimlaan 1, B-2020 Antwerp, Belgium
	}

	\author{%
		Puya Latafat\thanks{%
			\TheAddressKU\ {\tt
				\{%
					\href{mailto:puya.latafat@esat.kuleuven.be}{puya.latafat},%
					\href{mailto:panos.patrinos@esat.kuleuven.be}{panos.patrinos}%
				\}%
				\href{mailto:puya.latafat@esat.kuleuven.be,panos.patrinos@esat.kuleuven.be}{@esat.kuleuven.be}%
			}%
		}\and
		Andreas Themelis\thanks{%
			\TheAddressKUJ\ {\tt
				\href{mailto:andreas.themelis@ees.kyushu-u.ac.jp}{andreas.themelis@ees.kyushu-u.ac.jp}%
			}%
		}\and
		Masoud Ahookhosh\thanks{%
			\TheAddressUA\ {\tt
				\href{mailto:masoud.ahookhosh@uantwerp.be}{masoud.ahookhosh@uantwerp.be}%
			}%
		}\and
		Panagiotis Patrinos\footnotemark[2]%
	}%

	\title{\texorpdfstring{\TheTitle\thanks{\TheFunding}}{\TheShortTitle}}
	\date{}

\begin{document}
	
	\maketitle
	\begin{abstract}%

		We introduce two algorithms for nonconvex regularized finite sum minimization, where typical Lipschitz differentiability assumptions are relaxed to the notion of relative smoothness.
		The first one is a Bregman extension of Finito/MISO \cite{defazio2014finito,mairal2015incremental}, studied for fully nonconvex problems when the sampling is \bl{randomized}, or under convexity of the nonsmooth term when it is essentially cyclic.
		The second algorithm is a low-memory variant, in the spirit of SVRG \cite{johnson2013accelerating} and SARAH \cite{nguyen2017sarah}, that also allows for fully nonconvex formulations.
		Our analysis is made remarkably simple by employing a Bregman Moreau envelope as Lyapunov function.
		In the randomized case, linear convergence is established when the cost function is strongly convex, yet with no convexity requirements on the individual functions in the sum.
		For the essentially cyclic and low-memory variants, global and linear convergence results are established when the cost function satisfies the Kurdyka-\L ojasiewicz property.
	\end{abstract}%
	
	\begin{keywords}\TheKeywords \end{keywords}%
	\begin{AMS}\TheAMSsubj \end{AMS}%

	\section{Introduction}%

		We consider the following regularized finite sum minimization
		\begin{equation*}\tag{P}\label{eq:P}
			\textstyle
			\minimize_{x\in\R^n}{
				\varphi(x)
			{}\coloneqq{}
				\tfrac1N\sum_{i=1}^N f_i(x)
				{}+{}
				g(x)
			}\quad
			\stt x\in\overline C,
		\end{equation*}
		where \(\overline C\) denotes the closure of \(C\coloneqq\bigcap_{i=1}^N\interior\dom h_i\), for some convex functions \(h_i\), \(i\in[N]\coloneqq\set{1,\dots,N}\).
		Our goal in this paper is to study such problems without imposing convexity assumptions on $f_i$ and $g$, and in a setting where $f_i$ are differentiable but their gradients need not be Lipschitz continuous.
		Our full setting is formalized in \cref{ass:basic}.
		
		To relax the Lipschitz differentiability assumption, we adopt the notion of smoothness relative to a distance-generating function \cite{bauschke2017descent}, and following \cite{lu2018relatively} we will use the terminology of \emph{relative smoothness}.
		Despite the lack of Lipschitz differentiability, in many applications the involved functions satisfy a descent property where the usual quadratic upper bound is replaced by a Bregman distance (cf. \cref{thm:LipIneq,defin:Bregman}).
		Owing to this property, Bregman extensions for many classical schemes have been proposed \cite{bauschke2017descent,lu2018relatively,bauschke2019linear,teboulle2018simplified,ochs2019nonsmooth,ahookhosh2021bregman}.
		
		In the setting of finite sum minimization, the incremental aggregated algorithm PLIAG was proposed recently \cite{zhang2021proximal} as a Bregman variant of the incremental aggregated gradient method \cite{bertsekas2000gradient,blatt2007convergent,vanli2018global}.
		The analysis of PLIAG is limited to the convex case and requires restrictive assumptions for the Bregman kernel \cite[\bl{Thm. 1, Assump. 8}]{zhang2021proximal}.
		Stochastic  mirror descent (SMD) is another relevant algorithm which can tackle more general stochastic optimization problems.
		SMD may be viewed as a Bregman extension of the stochastic (sub)gradient method and has long been studied \cite{nemirovski2009robust,tseng2008accelerated,beck2003mirror,nedic2014stochastic}.
		More recently, \cite{hanzely2018fastest} studied SMD for convex and relatively smooth formulations, and (sub)gradient versions have been analyzed under relative continuity in a convex setting \cite{lu2019relative}, as well as relative weak convexity \cite{zhang2018convergence,davis2018stochastic}.
		
		Motivated by these recent works, we propose a Bregman extension of the popular Finito/\allowbreak MISO algorithm \cite{defazio2014finito,mairal2015incremental} in a fully nonconvex setting and with very general sampling strategies that will be made percise shortly after.
		In a nutshell, our analysis revolves around the fact that, regardless of the index selection strategy, the function \(\func{\LL}{\R^n\times\R^{nN}}{\Rinf}\) defined as
		\begin{equation}\label{eq:L}
			\LL(z,\bm s)
		{}\coloneqq{}
		\textstyle	\varphi(z)+\sum_{i=1}^N\D**_i\bigl(s_i,\nabla\h*_i(z)\bigr),
		\end{equation}
		where \(\conj{\h*_i}\) denotes the convex conjugate of \(\h*_i\coloneqq\nicefrac{\h_i}{\gamma_i}-\nicefrac{f_i}{N}\), monotonically decreases along the iterates \(\seq{z^k,\bm s^k}\) generated by \cref{alg:Finito} (see \cref{ass:basic} for the requirements on $\h_i,f_i$).
		Our methodology leverages an interpretation of Finito/MISO as a block-coordinate algorithm that was observed in \cite{latafat2021block} in the Euclidean setting.
		In fact, the analysis is here further simplified after noticing that the smooth function can be ``hidden'' in the distance-generating function, resulting in a Lyapunov function \(\LL\) that can be expressed as a \emph{Bregman Moreau envelope} (cf. \cref{thm:equivBC}).
		
		\begin{algorithm}[tb]
			\caption{Bregman Finito/MISO (BFinito) for the regularized finite sum minimization \eqref{eq:P}}%
			\label{alg:Finito}%
			\begin{algorithmic}[1]
			\Require
				\begin{tabular}[t]{@{}l@{}}
					Legendre kernels \(\h_i\) such that \(f_i\) is \(L_{f_i}\)-smooth relative to \(\h_i\)
				\\
					stepsizes \(\gamma_i\in(0,\nicefrac{N}{L_{f_i}})\)
				\\
					initial point \(x^{\rm init}\in C\coloneqq\bigcap_{i=1}^N\interior\dom h_i\)
				\end{tabular}
			\Initialize
				\begin{tabular}[t]{@{}l@{}}
					table \(\bm s^0 =(s_1^0,\ldots,s_N^0)\in\R^{nN}\) of vectors \(s_i^0=\tfrac{1}{\gamma_i}\nabla\h_i(x^{\rm init})-\tfrac1N\nabla f_i(x^{\rm init})\)
				\\
					\(\R^n\)-vector \(\s^0=\sum_{i=1}^Ns_i^0\)
				\end{tabular}
			\item[\algfont{Repeat} for \(k=0,1,\ldots\) until convergence]
			\State\label{state:Finito:z}%
				Compute
				\(
					z^k
				{}\in{}
					\argmin_{w\in\R^n}\set{
						g(w)
						{}+{}
						\sum_{i=1}^N\tfrac{1}{\gamma_i}\h_i(w)
						{}-{}
						\innprod*{\s^k}{w}
					}
				\)
			\State\label{state:Finito:I}
				Select a subset of indices \(\mathcal I^{k+1}\subseteq[N]\coloneqq\{1,\ldots,N\}\) and update the table \(\bm s^{k+1}\) as follows:
				\[
					s_i^{k+1}
				{}={}
					\begin{ifcases}
						\tfrac{1}{\gamma_i}\nabla \h_i(z^k)-\tfrac1N\nabla f_i(z^k) & i\in\mathcal I^{k+1}
					\\[3pt]
						s_i^k\otherwise
					\end{ifcases}
				\]
			\State\label{state:Finito:s+}%
				Update the vector
				\(
					\s^{k+1}
				{}={}
					\s^k
					{}+{}
					\sum_{i\in\mathcal I^{k+1}}(s_i^{k+1}-s_i^k)
				\)
			\item[\algfont{Return}]
				\(z^k\)
			\end{algorithmic}
		\end{algorithm}
		
		We cover a wide range of sampling strategies for the index set \(\mathcal I^{k+1}\) at \cref{state:Finito:I}, which we can summarize into the following two umbrella categories:
		\begin{align*}
		\tag{\ensuremath{\mathcal S_1}}\label{eq:RS}
			\text{\sc \bl{Randomized} rule: }
			&\quad
			\exists p_1,\ldots,p_N>0:
			~~
			\P{i\in\mathcal I^{k+1}}=p_i
			~~
			\forall k\in\N,i\in[N].
		\\
		\tag{\ensuremath{\mathcal S_2}}\label{eq:ECS}
			\text{\sc Essentially cyclic rule: }
			&\quad
			\textstyle
			\exists T>0:
			~~
			\bigcup_{t=1}^T\mathcal I^{k+t}=[N]
			~~
			\forall k\in\N.
		\end{align*}
		The randomized setting \eqref{eq:RS}, in which \(\P{}\) denotes the probability conditional to the knowledge at iteration \(k\), covers, for instance, a mini-batch strategy of size $b$.
		Another notable case is when each index $i$ is selected at random with probability $p_i$ independently of other indices.
		
		{%
		\setlength\belowdisplayshortskip{0pt}%
		\setlength\belowdisplayskip{0pt}%
			The essentially cyclic rule \eqref{eq:ECS} is also very general and has been considered by many authors \cite{tseng1987relaxation,tseng2001convergence,hong2017iteration,chow2017cyclic,xu2017globally}.
			Two notable special cases of single index selection rules complying with \eqref{eq:ECS} are the cyclic and shuffled cyclic sampling strategies:
			\begin{align*}
			\tag{\ensuremath{\mathcal S_2^{\text{\sc shuf}}}}\label{eq:shuffled}
				\text{\sc Shuffled cyclic rule:}
			&\quad
				\mathcal I^{k+1}
			{}={}
				\set{\pi_{\lfloor\nicefrac kN\rfloor}\bigl(\mod(k,N)+1\bigr)},
			\shortintertext{%
				where \(\pi_0,\pi_1,\ldots\) are permutations of the set of indices \([N]\) (chosen randomly or deterministically).
				When \(\pi_{\lfloor\nicefrac kN\rfloor}=\id\) one recovers the (plain) cyclic sampling rule
			}
			\tag{\ensuremath{\mathcal S_2^{\text{\sc cycl}}}}\label{eq:cyclic}
				\text{\sc Cyclic rule: }
			&\quad
				{\mathcal I}^{k+1}=\set{\mod(k,N)+1}.
			\end{align*}
			We remark that, in the cyclic case, our algorithm generalizes DIAG \cite{mokhtari2018surpassing} for smooth strongly convex problems, which itself may be seen as a cyclic variant of Finito/MISO.
		}%

		\subsection{Low-memory variant}%

			\begin{algorithm}[tb]%
				\caption{Low-memory Bregman Finito/MISO}%
				\label{alg:LM}%
				\begin{algorithmic}[1]
				\Require
					\begin{tabular}[t]{@{}l@{}}
						Legendre kernels \(\h_i\) such that \(f_i\) is \(L_{f_i}\)-smooth relative to \(\h_i\)
					\\
						stepsizes \(\gamma_i\in(0,\nicefrac{N}{L_{f_i}})\)
					\\
						initial point \(x^{\rm init}\in C\coloneqq\bigcap_{i=1}^N\interior\dom h_i\)
					\end{tabular}
				\Initialize
					\begin{tabular}[t]{@{}l@{}}
						\(\R^n\)-vector
						\(
							\tilde s^0
						{}={}
							\sum_{i=1}^N\bigl[
								\tfrac{1}{\gamma_i}\nabla\h_i(x^{\rm init})-\frac{1}{N}\nabla f_i(x^{\rm init})
							\bigr]
						\)
					\\
						set of selectable indices \(\mathcal K^0=\emptyset\)
					\Comment{Conventionally set to \(\emptyset\) so as to start with a full update}%
					\end{tabular}
				\item[\algfont{Repeat for \(k=0,1,\ldots\)} until convergence]
				\State\label{state:LM:z}%
					\(
						z^k
					{}\in{}
						\argmin_{w\in\R^n}\set{
							g(w)
							{}+{}
							\sum_{i=1}^N\tfrac{1}{\gamma_i}\h_i(w)
							{}-{}
							\innprod*{\s^k}{w}
						}
					\)
				\If{~\(\mathcal K^k=\emptyset\)~}\label{state:LM:full}%
					\Comment{No index left to be sampled: full update}
					\renewcommand{\ralignToWidthOf}{\ensuremath{\mathcal I^{k+1}=\mathcal K^{k+1}=[N]\quad}}%
					\State\label{state:LM:IK+_full}%
						\(
							\ralign{\mathcal I^{k+1}=\mathcal K^{k+1}=[N]}
						\)
					\Comment{activate all indices and reset the selectable indices}
					\State\label{state:LM:tz_full}%
						\(
							\ralign{\tilde z^k=z^k}
						\)
					\Comment{store the full update \(z^k\)}
					\State\label{state:LM:ts+_full}%
						\(
							\tilde s^{k+1}
						{}={}
							\sum_{i=1}^N\bigl[
								\tfrac{1}{\gamma_i}\nabla\h_i(z^k)-\tfrac1N\sum_{i=1}^{N}\nabla f_i(z^k)
							\bigr]
						\)
				\Else
					\renewcommand{\ralignToWidthOf}{\text{select a nonempty subset of indices \(\mathcal I^{k+1}\subseteq\mathcal K^k\)}}%
					\State\label{state:LM:I+}%
						\ralign{select a nonempty subset of indices \(\mathcal I^{k+1}\subseteq\mathcal K^k\)}%
					\Comment{select among the indices not yet sampled}%
					\State\label{state:LM:K+}%
						\ralign{\(\mathcal K^{k+1}=\mathcal K^k\setminus\mathcal I^{k+1}\)}%
					\Comment{update the set of selectable indices}
					\State\label{state:LM:tz}%
						\(
							\tilde z^k=\tilde z^{k-1}
						\)
					\State\label{state:LM:ts+}%
					\(
						\s^{k+1}
					{}={}
						\s^k
						{}+{}
						\sum_{i\in\mathcal I^{k+1}}\left[
							\bigl(
								\tfrac{1}{\gamma_i}
								\nabla\h_i(z^k)
								{}-{}
								\tfrac1N
								\nabla f_i(z^k)
							\bigr)
							{}-{}
							\bigl(
								\tfrac{1}{\gamma_i}
								\nabla\h_i(\tilde z^k)
								{}-{}
								\tfrac1N
								\nabla f_i(\tilde z^k)
							\bigr)
						\right]
					\)
				\EndIf
				\item[\algfont{Return}]
					\(\tilde z^k\)
				\end{algorithmic}
			\end{algorithm}
			
			One iteration of \cref{alg:Finito} involves the computation of \(z^k\) at \cref{state:Finito:z} and that of the gradients \(\nabla(\nicefrac{\h_i}{\gamma_i}-\nicefrac{f_i}{N})\), \(i\in\mathcal I^{k+1}\), at \cref{state:Finito:s+}.
			Consequently, the overall complexity of each iteration is independent of the number \(N\) of functions appearing in problem \eqref{eq:P}, and is instead proportional to the number of sampled indices, which the user is allowed to upper bound to any integer between 1 and \(N\).
			As is the case for all incremental gradient methods, the low iteration cost comes at the price of having to store in memory a table \(\bm s^k\) of \(N\) many \(\R^n\) vectors, which can become problematic when \(N\) grows large.
			Other incremental algorithms for convex optimization such as IAG \cite{bertsekas2000gradient,blatt2007convergent,vanli2018global}, IUG \cite{tseng2014incrementally}, SAG \cite{schmidt2017minimizing}, and SAGA \cite{defazio2014saga}, can considerably reduce memory allocation from \(O(nN)\) to \(O(n)\) in applications such as logistic regression and lasso where the gradients $\nabla f_i$ can be expressed as scaled versions of the data vectors.
			Despite the favorable performance of the Finito/MISO algorithm on such problems as observed in \cite{defazio2014saga}, this memory reduction trick can not be employed due to the fact that the vectors $s_i$ stored in the table depend not only on the gradients, but also on the vectors \(\nabla \h_i(z^k)\).
			Nevertheless, inspired by the popular stochastic methods SVRG \cite{johnson2013accelerating,xiao2014proximal} and SARAH \cite{nguyen2017sarah}, by suitably interleaving incremental and full gradient evaluations it is possible to completely waive the need of a memory table and match the \(O(n)\) storage requirement.
			
			In a nutshell, after a full update --- which in \cref{alg:Finito} corresponds to selecting $\mathcal I^{k+1}=[N]$ --- all vectors \(s_i^{k+1}\) in the table only depend on variable \(z^k\) computed at \cref{state:Finito:z}, until \(i\) is sampled again.
			As long as full gradient updates are frequent enough so that no index is sampled twice in between, it thus suffices to keep track of \(z^k\in\R^n\) instead of the table \(\bm s^k\in\R^{nN}\).
			The variant is detailed in \cref{alg:LM}, in which \(\mathcal K^k\subseteq[N]\) keeps track of the indices that have not yet been sampled between full gradient updates (and is thus reset whenever such full steps occur, cf. \cref{state:LM:IK+_full}).
			Vector \(\tilde z^k\in\R^n\) is equal to \(z^k\) corresponding to the latest full gradient update (cf. \cref{state:LM:tz_full}) and acts as a low-memory surrogate of the table \(\bm s^k\) of \cref{alg:Finito}.
			Similarly to SVRG and SARAH, this reduction in the storage requirements comes at the cost of an extra gradient evaluation per sampled index (cf. \cref{state:LM:ts+}).
			
			Since full gradient updates correspond to selecting all indices, \cref{alg:LM} may be viewed as \cref{alg:Finito} with an essentially cyclic sampling rule of period \(N\), a claim that will be formally shown in \cref{thm:equivLM}.
			In fact, not only does it naturally inherit all the convergence results, but its particular sampling strategy also allows us to waive convexity requirements on \(g\) that are necessary for more general essentially cyclic rules.
		
		\subsection{Contributions}%
			As a means to informally summarize the content of the paper, in \cref{table:summary} we synopsize the convergence results of the two algorithms.
			\begin{enumerate}[%
				leftmargin=0pt,
				label={\arabic*.},
				align=left,
				labelwidth=1em,
				itemindent=\labelwidth+\labelsep,
			]
			\item
				To the best of our knowledge, this is the first analysis of an incremental aggregated method in a fully nonconvex setting and without Lipschitz differentiability assumptions.
				Our analysis, surprisingly simple and yet \bl{covering} randomized and essentially cyclic samplings altogether, relies on a sure descent property on the Bregman Moreau envelope (cf. \cref{thm:sure}).
			\item
				We propose a novel low-memory variant of the (Bregman) Finito/MISO algorithm, that in the spirit of SVRG \cite{johnson2013accelerating,xiao2014proximal} and SARAH \cite{nguyen2017sarah} alternates between incremental steps and a full proximal gradient step.
				It is highly interesting even in the Euclidean case, as it can accommodate fully nonconvex formulations while maintaining a $O(n)$ memory requirement.
			\item
				Linear convergence of \cref{alg:Finito} in the randomized case is established when the cost function \(\varphi\) is strongly convex, yet with no convexity requirement on \(f_i\) or \(g\).
				To the best of our knownledge, this is a novelty even in the Euclidean case, for all available results are bound to strong convexity of each term $f_i$ in the sum; see \eg \cite{defazio2014finito,mairal2015incremental,mokhtari2018surpassing,latafat2021block,qian2019miso}.
				\bl{%
					This type of assumption has also been considered in \cite{allenzhu2016improved} for the case of SVRG.
					Although in practice it is hardly ever the case that \(f\) is strongly convex without each \(f_i\) also being (strongly) convex, our analysis being agnostic to the decomposition leads to tighter and more general results.%
				}%
			\item
				We leverage the Kurdyka-\L ojasiewicz (KL) property to establish global (as apposed to subsequential) convergence as well as linear convergence, for \cref{alg:Finito} with (essentially) cyclic sampling and for the low-memory \cref{alg:LM}.
			\end{enumerate}
			
			\begin{table}[tb]
				\makeatletter
					\newcommand{\KL}{{\sc kl}\@ifstar{\(_{\nicefrac12}\)}{}}
					\renewcommand{\E}{\mathbb E}%
				\makeatother
				\centering

					\footnotesize
					\renewcommand{\arraystretch}{1.2}%
					\setlength\tabcolsep{5pt}%
					\newcommand{\header}[1]{\mc{1}{@{}c@{}}{\bf #1}}%
					\newcommand{\EITHER}{either~~\ignorespaces}%
					\newcommand{\OR}{\fillwidthof[r]{either}{or}~~\ignorespaces}%
					\let\mc\multicolumn%
					\newcommand{\mr}[2][2]{\multirow{#1}{*}{\begin{tabular}{@{}c@{}}#2\end{tabular}}}%
					\newcommand{\Mr}[2][3]{\mr[#1]{#2}}%
					\renewcommand{\phi}{\varphi}%
					\newcommand{\cost}{\ensuremath{\varphi+\indicator_{\overline C}}}%
					\newcommand{\club}{\bl{\(^{\clubsuit}\)}}%
				\begin{tabular}{@{}| c | c | l|c|l |@{}}
					\header{}	& \header{\clap{Sampling}}	& \header{Requirements ~\normalfont (additionally to \cref{ass:basic})}	& \header{\clap{Property}}	& \header{Reference}	\\
				\hline
				%
				%
					\Mr{$z^k$\\bounded}
						& \Mr{any}			& \Mr{\cost\ level bounded}														& \Mr{sure}	& \Mr{\cref{thm:bounded}}	\\
						& 					&																				&			&							\\
						&					&																				&			&							\\
				\hline\hline
				%
				%
					\Mr{$\phi(z^k)$\\convergent}
						& \ref{eq:RS}		& 																				& a.s.		& \cref{thm:RS:cost}	\\\cline{2-5}
						& \ref{eq:ECS}		&$C=\R^n$;~ $g$ cvx;~ $\h_i$ loc str cvx smooth;~ $\varphi$ level bounded\club	& \mr{sure}	& \cref{thm:ECS:subseq}	\\\cline{2-3}\cline{5-5}
						& LM				& 																				&			& \cref{thm:LM:subseq}	\\
				\hline\hline
				%
				%
					\Mr[4]{$\omega(z^k)$\\stationary}
						& \mr{\ref{eq:RS}}	& \EITHER 	$C=\R^n$															& \mr{a.s.}	& \cref{thm:RS:stationary}			\\\cline{5-5}
						&					& \OR 		$\dom\h_i$ closed; $\varphi$ cvx									& 			& \cref{thm:RS:convex:stationary}	\\\cline{2-5}
						& \ref{eq:ECS}		& $C=\R^n$;~ $g$ cvx;~ $\h_i$ loc str cvx smooth;~ $\varphi$ level bounded\club	& \mr{sure}	& \cref{thm:ECS:subseq}				\\\cline{2-3}\cline{5-5}
						& LM				& $C=\R^n$																		&			& \cref{thm:LM:subseq}				\\
				\hline\hline
				%
				%
					\Mr[4]{$z^k$\\ convergent}
						& \mr{\ref{eq:RS}}	& \EITHER 	$C=\R^n$;~ $\phi$ cvx												& \mr{a.s.}	& \mr{\cref{thm:RS:convex:seq}}	\\
						&					& \OR 		\cref{ass:d};~ \cost\ cvx level bounded								&  			& 								\\\cline{2-5}
						& \ref{eq:ECS}		& \cref{ass:KL};~ $g$ cvx														& \mr{sure}	& \cref{thm:ECS:global:z}		\\\cline{2-3}\cline{5-5}
						& LM				& \cref{ass:KL}																	&			& \cref{thm:LM:global:z}		\\
				\hline\hline
				%
				%
					\Mr[3]{$\phi(z^k)$ and $z^k$\\linearly\\convergent}
						& \ref{eq:RS}		& $C=\R^n$;~ $\varphi$ str cvx; $\h_i$ loc smooth								& $\E$		& \cref{thm:RS:strconvex}		\\\cline{2-5}
						& \ref{eq:ECS}		& \cref{ass:KL};~ $\phi$ \KL*;~ $g$ cvx											& \mr{sure}	& \cref{thm:ECS:global:linear}	\\\cline{2-3}\cline{5-5}
						& LM				& \cref{ass:KL};~ $\phi$ \KL* 													&			& \cref{thm:LM:global:linear}	\\
				\hline
				\multicolumn{5}{@{}l@{}}{%
					\bl{\footnotesize
						\(^\clubsuit\)%
						Level boundedness is not necessary if \(h_i\) are \emph{globally} smooth and strongly convex (cf. \cref{thm:ECS:subseq}).%
					}%
				}%
				\end{tabular}%
				\caption[Summary of the convergence results.]{%
					Summary of the convergence results for \cref{alg:Finito} with
					randomized \bl{rule} \eqref{eq:RS} and
					essentially cyclic \bl{rule} \eqref{eq:ECS},
					and for the low-memory variant of \cref{alg:LM} {\em (LM)}.
					Claims are either
					{\em sure},
					almost sure {\em (a.s.)},
					or in expectation {\em (\(\E\))}.
					\\
					{\sc Other abbreviations:}~
					{\em loc}: locally;~
					{\em cvx}: convex;~
					{\em str}: strongly;~
					{\em smooth}: Lipschitz differentiable;~
					\(\omega\): set of limit points;~
					\bl{\KL*}: Kurdyka-\L ojasiewicz property with exponent \bl{\(\nicefrac12\)}.%
				}%
				\label{table:summary}%
			\end{table}

		\subsection{Organization}%
			We conclude this section by introducing some notational conventions.
			The problem setting is formally described in \cref{sec:preliminaries} together with a list of related definitions and known facts involving Bregman distances, relative smoothness, and proximal mapping.
			\Cref{sec:BC} offers an alternative interpretation of \cref{alg:Finito} as the block-coordinate Bregman proximal point \cref{alg:BCPP}, which majorly simplifies the analysis, addressed in \cref{sec:conv}.
			Some auxiliary results are deferred to 
			\cref{sec:appendix:aux}.
			\Cref{sec:simulations} applies the proposed algorithms to sparse phase retrieval problems, and \cref{sec:conclusions} concludes the paper.

		\subsection{Notation}%
			The set of real and extended-real numbers are \(\R\coloneqq(-\infty,\infty)\) and \(\Rinf\coloneqq\R\cup\set\infty\), while the positive and strictly positive reals are \(\R_+\coloneqq[0,\infty)\) and \(\R_{++}\coloneqq(0,\infty)\).
			With \(\id\) we indicate the identity function \(x\mapsto x\) defined on a suitable space.
			We denote by \(\innprod{{}\cdot{}}{{}\cdot{}}\) and \(\|{}\cdot{}\|\) the standard Euclidean inner product and the induced norm.
			For a vector $\bm w=(w_1,\ldots,w_r)\in\R^{\sum_i n_i}$, $w_i\in\R^{n_i}$ is used to denote its $i$-th block coordinate.
			\(\interior E\) and \(\boundary E\) respectively denote the interior and boundary of a set \(E\), and for a sequence \(\seq{x^k}\) we write \(\seq{x^k}\subseteq E\) to indicate that \(x^k\in E\) for all \(k\in\N\).
			We say that \(\seq{x^k}\) converges at \DEF{\(Q\)-linear rate} (resp. \DEF{\(R\)-linear rate}) to a point \(x\) if there exists \(c\in(0,1)\) such that \(\|x^{k+1}-x\|\leq c\|x^k-x\|\) (resp. \(\|x^k-x\|\leq\rho c^k\) for some \(\rho>0\)) holds for all \(k\in\N\).
			
			
			We use the notation \(\ffunc H{\R^n}{\R^m}\) to indicate a mapping from each point $x\in\R^n$ to a subset $H(x)$ of $\R^m$.
			The \DEF{graph} of \(H\) is the set
			\(
				\graph H
			{}\coloneqq{}
				\set{(x,y)\in\R^n\times\R^m}[
					y\in H(x)
				]
			\).
			We say that \(H\) is \DEF{outer semicontinuous (osc)} if \(\graph H\) is a closed subset of \(\R^n\times\R^m\), and \DEF{locally bounded} if for every bounded \(U\subset\R^n\) the set \(\bigcup_{x\in U}H(x)\) is bounded.
			
			The \DEF{domain} and \DEF{epigraph} of an extended-real-valued function \(\func{h}{\R^n}{\Rinf}\) are the sets
			\(
				\dom h
			{}\coloneqq{}
				\set{x\in\R^n}[
					h(x)<\infty
				]
			\)
			and
			\(
				\epi h
			{}\coloneqq{}
				\set{(x,\alpha)\in\R^n\times\R}[
					h(x)\leq\alpha
				]
			\). Function
			\(h\) is said to be \DEF{proper} if \(\dom h\neq\emptyset\), and \DEF{lower semicontinuous (lsc)} if \(\epi h\) is a closed subset of \(\R^{n+1}\).
			We say that \(h\) is \emph{level bounded} if its \(\alpha\)-sublevel set
			\(
				\lev_{\leq\alpha}h
			{}\coloneqq{}
				\set{x\in\R^n}[
					h(x)\leq\alpha
				]
			\) is bounded for all \(\alpha\in\R\).
			The \DEF{conjugate} of \(h\), is defined by
			\(
				\conj h(y)
			{}\coloneqq{}
				\sup_{x\in\R^n}\set{\innprod yx-h(x)}
			\).
			The indicator function of a set $E\subseteq\R^n$ is denoted by \(\indicator_E\), namely \(\indicator_E(x)=0\) if \(x\in E\) and \(\infty\) otherwise.
			
			We denote by \(\ffunc{\hat\partial h}{\R^n}{\R^n}\) the \DEF{regular subdifferential} of \(h\), where
			\begin{equation*}
				v\in\hat\partial h(\bar x)
			\quad\Leftrightarrow\quad
				\liminf_{\bar x\neq x\to\bar x}{
					\frac{h(x)-h(\bar x)-\innprod{v}{x-\bar x}}{\|x-\bar x\|}
				}
			{}\geq{}
				0.
			\end{equation*}
			A necessary condition for local minimality of \(x\) for \(h\) is \(0\in\hat\partial h(x)\), see \cite[Th. 10.1]{rockafellar2009variational}.
			The (limiting) \DEF{subdifferential} of \(h\) is \(\ffunc{\partial h}{\R^n}{\R^n}\), where
			\(
				v\in\partial h(x)
			\)
			iff \(x\in\dom h\) and there exists a sequence \(\seq{x^k,v^k}\subseteq\graph\hat\partial h\) such that
			\(
				(x^k,h(x^k),v^k)
			{}\to{}
				(x,h(x),v)
			\)
			as \(k\to\infty\).
			Finally, the set of $r$ times continuously differentiable functions from $X$ to $\R$ is denoted by $\C^r(X)$.

\renewcommand{\BregmanKernel}{h}%
	\section{Problem setting and preliminaries}\label{sec:preliminaries}%

		Throughout this paper, problem \eqref{eq:P} is studied under the following assumptions.
		
		\begin{assumption}[basic requirements]\label{ass:basic}
			In problem \eqref{eq:P},
			\begin{enumeratass}
			\item\label{ass:f}%
				\(\func{f_i}{\R^n}{\Rinf}\) are \(L_{f_i}\)-smooth relative to Legendre kernels \(\h_i\) (\cref{defin:LegendreKernel,defin:relsmooth});%
			\item\label{ass:g}%
				\(\func{g}{\R^n}{\Rinf}\) is proper and lower semicontinuous (lsc);
			\item\label{ass:argmin}%
				a solution exists: \(\argmin\set{\varphi(x)}[x\in\overline C]\neq\emptyset\);
			\item\label{ass:rangeT}%
				for given \(\gamma_i\in(0,\nicefrac{N}{L_{f_i}})\), \(i\in[N]\), it holds that for any \(s\in\R^n\)
				\begin{equation}\label{eq:innerprox}
					T(s)
				{}\coloneqq{}
					\argmin_{w\in\R^n}\set{\textstyle
						g(w)
						{}+{}
						\sum_{i=1}^N\tfrac{1}{\gamma_i}\h_i(w)
						{}-{}
						\innprod*{s}{w}
					}
				{}\subseteq{}
					C.
				\end{equation}
			\end{enumeratass}
		\end{assumption}
		
		As it will become clear in \cref{sec:BC}, the subproblem \eqref{eq:innerprox} is in fact a reformulation of a (Bregman) proximal mapping. \Cref{ass:rangeT} excludes boundary points from \(\range T\).
		This is a standard assumption that usually holds in practice \cite{bolte2018first,solodov2000inexact}, \eg when $g$ is convex or when the intersection of $\dom\h_i$, $i\in[N]$, is an open set.

		\begin{definition}[Bregman distance]\label{defin:Bregman}%
			For a convex function \(\func{\h}{\R^n}{\Rinf}\) that is continuously differentiable on \(\interior\dom\h\neq\emptyset\), the \DEF{Bregman distance} \(\func{\D}{\R^n\times\R^n}{\Rinf}\) is
		defined as
			\begin{equation}\label{eq:bregman}
				\D(x,y)
			{}\coloneqq{}
				\begin{ifcases}
					\h(x)-\h(y)-\innprod*{\nabla\h(y)}{x-y} & y\in\interior\dom\h
				\\
					\infty\otherwise.
				\end{ifcases}
			\end{equation}
			Function \(\h\) will be referred to as a \DEF{distance-generating function}.
		\end{definition}
		
		\begin{definition}[Legendre kernel]\label{defin:LegendreKernel}%
			A proper, lsc, and strictly convex function \(\func{\h}{\R^n}{\Rinf}\) with \(\interior\dom\h\neq\emptyset\) and such that \(\h\in\C^1(\interior\dom\h)\) is said to be a Legendre kernel if it is
			(i) \DEF{\(1\)-coercive}, \ie such that \(\lim_{\|x\|\to\infty}\nicefrac{\h(x)}{\|x\|}=\infty\), and
			(ii) \DEF{essentially smooth}, \ie if \(\|\nabla\h(x_k)\|\to\infty\) for every sequence \(\seq{x_k}\subseteq\interior\dom\h\) converging to a boundary point of \(\dom\h\).
		\end{definition}
		
		\begin{fact}
			\bl{Let \(\func{\h}{\R^n}{\Rinf}\) be a Legendre kernel, \(x\in\R^n\), and \(y,z\in\interior\dom\h\).
			Then:}%
			\begin{enumerate}
			\item\label{thm:inv1}%
				\(\conj\h\in\C^1(\R^n)\) is strictly convex and
				\(\nabla\h^{-1}=\nabla\conj\h\)
				\cite[Thm. 26.5 and Cor. 13.3.1]{rockafellar1970convex}.
			\item\label{thm:3P}%
				\(
					\D(x,z)=\D(x,y)+\D(y,z)+\innprod*{x-y}{\nabla\h(y)-\nabla\h(z)}
				\)
				\cite[Lem. 3.1]{chen1993convergence}.
			\item\label{thm:D*}%
				\(\D(y,z)=\BregmanD_{\conj{\h}}(\nabla\h(z),\nabla\h(y))\) \cite[Thm. 3.7(v)]{bauschke1997legendre}.
			\item\label{thm:Dhcoercive}%
				\(\D({}\cdot{},z)\) and \(\D(z,{}\cdot{})\) are level bounded \cite[Lem. 7.3(v)-(viii)]{bauschke2001essential}.
			\item\label{thm:Dboundary}%
				If \(\dom\h\) is closed and \(\D(x^k,y^k)\to0\) for some \(x^k\in\dom\h\) and \(y^k\in\interior\dom\h\), then \(\seq{x^k}\) converges to a point \(x\) iff so does \(\seq{y^k}\) \cite[Thm. 2.4]{solodov2000inexact}.
			\end{enumerate}
			Moreover, for any convex set \(\U\subseteq\interior\dom\h\) and \(u,v\in\U\) the following hold:
			\begin{enumerate}[resume]
			\item\label{thm:hstrconvex}%
				If \(\h\) is \(\mu_{\h,\U}\)-strongly convex on \(\U\), then
				\(
					\frac{\mu_{\h,\U}}{2}\|v-u\|^2
				{}\leq{}
					\D(v,u)
				{}\leq{}
					\frac{1}{2\mu_{\h,\U}}\|\nabla\h(v)-\nabla\h(u)\|^2
				\).%
			\item\label{thm:hC11}%
				If \(\nabla\h\) is \(\ell_{\h,\U}\)-Lipschitz on \(\U\), then
				\(
					\frac{1}{2\ell_{\h,\U}}\|\nabla\h(v)-\nabla\h(u)\|^2
				{}\leq{}
					\D(v,u)
				{}\leq{}
					\frac{\ell_{\h,\U}}{2}\|v-u\|^2
				\).%
			\end{enumerate}
		\end{fact}
		
		\begin{definition}[relative smoothness {\cite{bauschke2017descent}}]\label{defin:relsmooth}%
			We say that a proper, lsc function \(\func{f}{\R^n}{\Rinf}\) is \DEF{smooth relative to} a Legendre kernel \(\func{\h}{\R^n}{\Rinf}\) if \(\dom f\supseteq\dom\h\), and there exists \(L_f\geq0\) such that \(L_f\h\pm f\) are convex functions on \(\interior\dom\h\).
			We will simply say that \(f\) is \(L_f\)-smooth relative to $\h$ to make the modulus \(L_f\) explicit.
		\end{definition}
		
		\begin{fact}
			Let \(\func f{\R^n}{\Rinf}\) be \(L_f\)-smooth relative to a Legendre kernel \(\func\h{\R^n}{\Rinf}\).
			Then, \(f\in\C^1(\interior\dom\h)\) and the following hold:
			\begin{enumerate}
			\item\label{thm:LipIneq}%
				\(
					\bigl|f(y)-f(x)-\innprod*{\nabla f(x)}{y-x}\bigr|
				{}\leq{}
					L_f\D(y,x)
				\)
				for all \(x,y\in\interior\dom\h\).
			\item\label{thm:C2relsmooth}%
				\(
					-L_f\nabla^2\h
				{}\preceq{}
					\nabla^2 f
				{}\preceq{}
					L_f\nabla^2\h
				\)
				on \(\interior\dom\h\), provided that \(f,h\in\mathcal{C}^2(\interior\dom\h)\).
			\item\label{thm:fC11}
				If \(\nabla\h\) is Lipschitz continuous with modulus \(\ell_{\h,\U}\) on a convex set \(\U\), then so is \(\nabla f\) with modulus \(\ell_{f,\U}=L_f\ell_{\h,\U}\) \cite[Prop. 2.5(ii)]{ahookhosh2021bregman}.
			\end{enumerate}
		\end{fact}

		Relative to a Legendre kernel \(\func{\h}{\R^{n}}{\Rinf}\), the \DEF{Bregman proximal mapping} of \(\psi\) is the set-valued map \(\ffunc{\prox_\FBEG}{\interior\dom\h}{\R^{n}}\) given by
		\begin{align}\label{eq:bprox}
			\prox_\FBEG(\bl x)
		{}\coloneqq{} &
			\argmin_{\bl z\in\R^{n}}\set{
				\psi(\bl z)+\D(\bl z,\bl x)
			},
		\shortintertext{
			and the corresponding \DEF{Bregman-Moreau envelope} is \(\func{\FBEG}{\R^{n}}{[-\infty,\infty]}\) defined as
		}
		\label{eq:bdre}
			\FBEG(\bl x)
		{}\coloneqq{} &
			\inf_{\bl z\in\R^{n}}\set{\psi(\bl z)+\D(\bl z,\bl x)}.
		\end{align}
		
		\begin{fact}[regularity properties of \(\prox_\FBEG\) and \(\FBEG\) \cite{kan2012moreau}]\label{thm:facts}%
			The following hold for a Legendre kernel \(\func{h}{\R^n}{\Rinf}\) and a proper, lsc, lower bounded function \(\func{\psi}{\R^n}{\Rinf}\):
			\begin{enumerate}
			\item\label{thm:osc}%
				\(\prox_\FBEG\) is locally bounded, compact-valued, and outer semicontinuous on \(\interior\dom\h\).
			\item\label{thm:C0}%
				\(\FBEG\) is real-valued and continuous on \(\interior\dom\h\); in fact, it is locally Lipschitz if so is \(\nabla\h\).
			\end{enumerate}
		\end{fact}
		
		\begin{fact}[relation between \(\psi\) and \(\FBEG\)]\label{thm:Moreau:basic}%
			Let \(\h\) be a Legendre kernel and \(\func{\psi}{\R^n}{\Rinf}\) be proper, lsc, and lower bounded on \(\overline{\dom\h}\).
			Then, for every \(x\in\interior\dom\h\), \(y\in\dom h\), and \(\bar x\in\prox_\FBEG(x)\)%
			\begin{enumerate}
			\item\label{thm:leq}%
				\(
					\FBEG(x)
				{}\defeq{}
					\psi(\bar x)
					{}+{}
					\D(\bar x,x)
				{}\leq{}
					\psi(y)
					{}+{}
					\D(y,x)
				\),
				and in particular \(\FBEG(x)\leq\psi(x)\);
			\item\label{thm:leq_convex}%
				if \(\psi\) is convex, then
				\(
					\FBEG(x)
				{}\leq{}
					\psi(y)
					{}+{}
					\D(y,x)
					{}-{}
					\D(y,\bar x)
				\)
				\cite[Lem. 3.1]{teboulle2018simplified}.
			\end{enumerate}
			Moreover, if \(\range\prox_\FBEG\subseteq\interior\dom\h\), then the following also hold \cite[Prop.3.3]{ahookhosh2021bregman}:
			\begin{enumerate}[resume]
			\item\label{thm:inf}
				\(
					\inf_{\overline{\dom\h}}\psi
				{}\leq{}
					\inf_{\interior\dom\h}\psi
				{}={}
					\inf\FBEG
				\)
				and
				\(
					\argmin\FBEG
				{}={}
					\argmin_{\interior\dom\h}\psi
				\).
			\item\label{thm:lb}
				\(\psi+\indicator_{\overline{\dom\h}}\) is level bounded iff so is \(\FBEG\).
			\end{enumerate}
		\end{fact}

\renewcommand{\BregmanKernel}{H}%
	\section{A block-coordinate interpretation}\label{sec:BC}%

		By introducing \(N\) copies of \(x\), problem \eqref{eq:P} can equivalently be written as
		\begin{equation}\label{eq:bigP}
			\minimize_{\bm x=(x_1,\ldots,x_N)\in\R^{nN}}{
				\Phi(\bm x)
			{}={}
				\overbracket*{
					\textstyle
					\tfrac1N\sum_{i=1}^Nf_i(x_i)
				}^{F(\bm x)}
				{}+{}
				\overbracket*{
					\textstyle
					\tfrac1N\sum_{i=1}^Ng(x_i)
					{}+{}
					\indicator_\Delta(\bm x)
				}^{G(\bm x)}
			}
		\quad
			\stt \bm x\in \overline C\times\cdots\times\overline C,
		\end{equation}
		where
		\(
			\Delta
		{}\coloneqq{}
			\set{\bm x=(x_1,\ldots,x_N)\in\R^{nN}}[
				x_1=x_2=\cdots=x_N
			]
		\)
		is the consensus set.
		The equivalence between \eqref{eq:bigP} and the original problem \eqref{eq:P} is formally established in \cref{thm:equivP}.
		Note that \cref{ass:f} implies that \(F\) as in \eqref{eq:bigP} is smooth with respect to the Legendre kernel
		\begin{equation}\label{eq:H}
			\func{\h}{\R^{nN}}{\Rinf}
		\quad\text{defined as}\quad
		\textstyle
			\h(\bm x)
		{}={}
			\sum_{i=1}^N\h_i(x_i),
		\end{equation}
		making Bregman forward-backward iterations
		\(
			\bm x^+
		{}\in{}
			\argmin\set*{
				\innprod*{\nabla F(\bm x)}{{}\cdot{}}
				{}+{}
				G({}\cdot{})
				{}+{}
				\tfrac1\gamma\D({}\cdot{},\bm x)
			}
		\)
		for some stepsize \(\gamma>0\) a suitable option to address problem \eqref{eq:bigP}.
		In fact, it can be easily verified that \(L_F=\tfrac1N\max_{i=1\ldots N}L_{f_i}\) is a smoothness modulus of \(F\) relative to \(\h\), indicating that fixed point iterations \(\bm x\gets\bm x^+\) under \cref{ass:basic} converge (in some sense to be made precise) to a stationary point of the problem whenever \(\gamma\in(0,\nicefrac{1}{L_F})\).
		Notice that a higher degree of flexibility can be granted by considering an \(N\)-uple of individual stepsizes \(\Gamma=(\gamma_1,\ldots,\gamma_N)\), giving rise to the forward-backward operator \(\ffunc{\T}{\R^{nN}}{\R^{nN}}\) in the Bregman metric \((\bm z,\bm x)\mapsto\sum_{i=1}^N\tfrac{1}{\gamma_i}\D_i(z_i,x_i)\), namely
		\begin{equation}\label{eq:T}
			\T(\bm x)
		{}\coloneqq{}
			\argmin_{\bm z\in\R^{nN}}\set{
			\textstyle
				F(\bm x)+\innprod*{\nabla F(\bm x)}{\bm z-\bm x}
				{}+{}
				G(\bm z)
				{}+{}
				\sum_{i=1}^N\tfrac{1}{\gamma_i}\D_i(z_i,x_i)
			}.
		\end{equation}
		This intuition is validated in the next result, which asserts that whenever the stepsizes \(\gamma_i\) are selected as in \cref{alg:Finito} the operator \(\T\) coincides with a proximal mapping on a suitable Legendre kernel function \(\h*\).
		This observation leads to a much simpler analysis of \cref{alg:Finito}, which will be shown to be a block-coordinate variant of a Bregman proximal point method.
		
		\begin{lemma}\label{thm:Legendre}%
			Suppose that \cref{ass:f} holds and let \(\gamma_i\in(0,\nicefrac{N}{L_{f_i}})\) be selected as in \cref{alg:Finito}.
			Then, \(\h*_i\coloneqq\frac{1}{\gamma_i}\h_i-\tfrac1N f_i\) (with the convention \(\infty-\infty=\infty\)) is a Legendre kernel with \(\dom\h*_i=\dom\h_i\), \(i\in[N]\), and thus so is
			the function
			\begin{equation}\label{eq:hatH}
				\func{\h*}{\R^{nN}}{\Rinf}
			\quad\text{defined as}\quad
			\textstyle
				\h*(\bm x)
			{}={}
				\sum_{i=1}^N\h*_i(x_i).
			\end{equation}
			Moreover, for any \((\bm z,\bm x)\in\R^{nN}\times\R^{nN}\) it holds that
			\begin{equation}\label{eq:M=PPM}
			\textstyle
				\Phi(\bm z)+\D*(\bm z, \bm x)
			{}={}
				F(\bm x)+\innprod*{\nabla F(\bm x)}{\bm z-\bm x}
				{}+{}
				G(\bm z)
				{}+{}
				\sum_{i=1}^N\tfrac{1}{\gamma_i}\D_i(z_i,x_i),
			\end{equation}
			and in particular the forward-backward operator \eqref{eq:T} satisfies
			\begin{equation}\label{eq:Tprox}
				\T(\bm x)
			{}={}
				\prox_\Phi^{\h*}(\bm x).
			\end{equation}
			When \cref{ass:basic} is satisfied, then the following also hold:
			\begin{enumerate}
			\item\label{thm:Legendre:hatHh}%
				\(
					\D*(\bm z,\bm x)
				{}\geq{}
					\sum_{i=1}^N(\tfrac{1}{\gamma_i}-\tfrac{L_{f_i}}{N})\D_i(z_i,x_i)
				\).
			\item\label{thm:Legendre:rangeT}%
				\(
					\prox_\Phi^{\h*}(\bm x)
				{}={}
					\set{(z,\cdots,z)}[
						z\in T(\sum_{i=1}^N\nabla\h*_i(x_i))
					]
				\),
				with \(T\) as in \eqref{eq:innerprox}, is a nonempty and compact subset of \(C\times\cdots\times C\) for any \(\bm x\in\interior\dom\h_1\times\cdots\times\interior\dom\h_N\).
			\item\label{thm:OC}%
				If \(\bm z\in\prox_\Phi^{\h*}(\bm x)\), then
				\(
					\nabla\h*(\bm x)
					{}-{}
					\nabla\h*(\bm z)
				{}\in{}
					\hat\partial\Phi(\bm z)
				\);
				the converse also holds when \(\Phi\) is convex.
			\item\label{thm:Legendre:smooth}%
				If $\nabla \h_i$ is $\ell_{\h_i,\U_i}$-Lipschitz on a convex set \(\U_i\subseteq\interior\dom\h_i\), then $\nabla\h*_i$ is $\ell_{\h*_i,\U_i}$-Lipschitz on $\U_i$ with
				\(
					\ell_{\h*_i,\U_i}
				{}\leq{}
					\big(\tfrac{1}{\gamma_i}+\tfrac{L_{f_i}}{N}\big)\ell_{\h_i,\U_i}
				\).
				If, in addition, $f_i - \tfrac{\mu_{f_i,\U_i}}2\|\cdot\|^2$ is convex on $\U_i$ for some $\mu_{f_i,\U_i}\in\R$, then
				\(
					\ell_{\h*_i}
				{}\leq{}
					\tfrac{\ell_{\h_i,\U_i}}{\gamma_i}-\tfrac{\mu_{f_i,\U_i}}{N}
				\).
			\item\label{thm:Legendre:strconvex}%
				If $\h_i$ is $\mu_{h_i,\U_i}$-strongly convex on a convex set \(\U_i\subseteq\dom\h_i\), then $\h*_i$ is $\mu_{\h*_i,\U_i}$-strongly convex on $\U_i$ with
				\(
					\mu_{\h*_i,\U_i}
				{}\geq{}
					\big(\tfrac{1}{\gamma_i}-\tfrac{L_{f_i}}{N}\big)\mu_{\h_i,\U_i}
				\).
			\end{enumerate}
			\begin{proof}
				The claims on \(\h*_i\) are shown in \cite[Thm. 4.1]{ahookhosh2021bregman}, and \eqref{eq:M=PPM} and \eqref{eq:Tprox} then easily follow.
				\begin{proofitemize}
				\item\ref{thm:Legendre:hatHh}~
					This is an immediate consenquence of \cref{thm:LipIneq}.
				\item\ref{thm:Legendre:rangeT}~
					Let \(\bm x\) be as in the statement, and observe that \(\bm x\in\interior\dom\h*\); nonemptyness and compactness of \(\prox_\Phi^{\h*}\) then follows from \cref{thm:osc}.
					Let now \(\bm u\in\prox_\Phi^{\h*}(\bm x)\) be fixed, and note that the consensus constraint encoded in \(\Phi\) ensures that \(u_i=u_j\) for all \(i,j\in[N]\).
					Thus,
					\begin{align*}
						u_i
					{}={} &
						\argmin_{w\in\R^n}\set{
							\Phi(w,\ldots,w)
							{}+{}
							\h*(w,\ldots,w)
							{}-{}
							\innprod*{\nabla\h*(\bm x)}{(w,\ldots,w)}
						}
					\\
					{}={} &
						\argmin_{w\in\R^n}\set{
							\textstyle
							\tfrac1N\sum_{i=1}^Nf_i(w)
							{}+{}
							g(w)
							{}+{}
							\bl{%
								\sum_{i=1}^N
								\bigl(
									\h*_i(w)
									{}-{}
									\innprod*{\nabla \h*_i(x_i)}{w}
								\bigr)
							}
						}
					\\
					{}={} &
						\argmin_{w\in\R^n\vphantom{_X}}\set{
							\textstyle
							g(w)
							{}+{}
							\sum_{i=1}^N{
								\tfrac{1}{\gamma_i}\h_i(w)
							}
							{}-{}
							\innprod*{
								\textstyle
								\sum_{i=1}^N\nabla\h*_i(x_i)
							}{w}
						}
					{}\overrel*{\eqref{eq:innerprox}}{}
						T\bigl(
							\textstyle
							\sum_{i=1}^N\nabla\h*_i(x_i)
						\bigr)
					{}\subseteq{}
						C,
					\end{align*}
					where the inclusion follows from \cref{ass:rangeT}.
				\item\ref{thm:OC}~
					Observe first that necessarily \(\bm x\in \interior\dom \h_i\times\cdots\times \interior\dom \h_N\), for otherwise no such \(\bm z\) exists.
					Moreover, from assertion \ref{thm:Legendre:rangeT} it follows that also \(\bm z\) belongs to such open set, onto which \(\h*\) is continuously differentiable.
					The claim then follows from the necessary condition for optimality of \(\bm z\) in the minimization problem \eqref{eq:bdre} --- which is also sufficient when \(\Phi\) is convex, for so is \(\Phi+\D*({}\cdot{},\bm x)\) in this case --- having
					\[
						0
					{}\in{}
						\hat\partial[\Phi+\D*({}\cdot{},\bm x)](\bm z)
					{}={}
						\hat\partial[\Phi+\h*-\innprod*{\nabla\h*(\bm x)}{{}\cdot{}}](\bm z)
					{}={}
						\hat\partial\Phi(\bm z)
						{}+{}
						\nabla\h*(\bm z)
						{}-{}
						\nabla\h*(\bm x).
					\]
					The last equality follows from \cite[Ex. 8.8(c)]{rockafellar2009variational}, owing to smoothness of \(\h*\) at \(\bm z\).
				\item\ref{thm:Legendre:smooth} and \ref{thm:Legendre:strconvex}~
					Observe that
					\begin{equation}\label{eq:hhat_LS}
					\let\frac\tfrac
						\frac{N-\gamma_i L_{f_i}}{N\gamma_i}\h_i
					{}\preceq{}
						\frac{N-\gamma_i L_{f_i}}{N\gamma_i}\h_i
						{}+{}
						\overbracket*{\frac1N(L_{f_i}\h_i-f_i)}^{\rm convex}
					{}={}
						\h*_i
					{}={}
						\frac{N+\gamma_i L_{f_i}}{N\gamma_i}\h_i
						{}-{}
						\overbracket*{\frac1N(L_{f_i}\h_i+f_i)}^{\rm convex}
					{}\preceq{}
						\frac{N+\gamma_i L_{f_i}}{N\gamma_i}\h_i,
					\end{equation}
					where for notational convenience we used the partial ordering ``\(\preceq\)'', defined as \(\alpha\preceq\beta\) iff \(\beta-\alpha\) is convex.
					The claimed moduli
					\(
						\ell_{\h*_i,\U_i}
					{}\leq{}
						\big(\frac{1}{\gamma_i}+\frac{L_{f_i}}{N}\big)\ell_{\h_i,\U_i}
					\)
					and
					\(
						\mu_{\h*_i,\U_i}
					{}\geq{}
						\big(\frac{1}{\gamma_i}-\frac{L_{f_i}}{N}\big)\mu_{\h_i,\U_i}
					\)
					are thus readily inferred.
					In case \(f_i\) is \(\mu_{f_i,\U_i}\)-strongly convex on \(\U_i\), we may write
					\[
						\h*_i
					{}={}
						\tfrac{\h_i}{\gamma_i}
						{}-{}
						\tfrac{\mu_{f_i,\U_i}}{2N}\|{}\cdot{}\|^2
						{}-{}
						\tfrac1N(\underbracket*{f_i-\tfrac{\mu_{f_i,\U_i}}{2}\|{}\cdot{}\|^2}_{\rm convex})
					{}\preceq{}
						\tfrac{\h_i}{\gamma_i}
						{}-{}
						\tfrac{\mu_{f_i,\U_i}}{2N}\|{}\cdot{}\|^2
					\]
					to obtain the tighter bound
					\(
						\ell_{\h*_i,\U_i}
					{}\leq{}
						\frac{\ell_{\h_i,\U_i}}{\gamma_i}-\frac{\mu_{f_i,\U_i}}{N}
					\).
				\qedhere
				\end{proofitemize}
			\end{proof}
		\end{lemma}
		
		\begin{algorithm}
			\caption{Block-coordinate proximal point formulation of \cref{alg:Finito}}%
			\label{alg:BCPP}%
			\begin{algorithmic}[1]
			\Require
				\begin{tabular}[t]{@{}l@{}}
					Legendre kernels \(\h_i\) such that \(f_i\) is \(L_{f_i}\)-smooth relative to \(\h_i\)
				\\
					stepsizes \(\gamma_i\in(0,\nicefrac{N}{L_{f_i}})\)
				\\
							initial point \(x^{\rm init}\in \cap_{i=1}^N\interior\dom h_i= C\)
				\end{tabular}
			\item[\algfont{\fillwidthof[l]{Initialize}{Denote}}]
				\(
					\bm x^0
				{}={}
					(x^{\rm init},\ldots,x^{\rm init})
				\),~~
				\(
					\h*_j
				{}\coloneqq{}
					\tfrac{1}{\gamma_i}\h_j-\tfrac{1}Nf_j
				\),~~
				\(
					\h*(\bm x)
				{}\coloneqq{}
					\sum_{i=1}^N\h*_i(x_i)
				\)
			\item[\algfont{Repeat} for \(k=0,1,\ldots\) until convergence]
			\State\label{state:BCPP:PP}%
				\(
					\bm u^k
				{}\in{}
					\argmin_{\bm w\in\R^{nN}}\set{
						\Phi(\bm w)
						{}+{}
						\h*(\bm w)
						{}-{}
						\innprod*{\nabla\h*(\bm x^k)}{\bm w}
					}
				{}={}
					\argmin_{\bm w\in\R^{nN}}\set{
						\Phi(\bm w)
						{}+{}
						\D*(\bm w,\bm x^k)
					}
				\)
			\State
				Select a subset of indices \(\mathcal J^{k+1}\subseteq[N]\)
			\State\label{state:BCPP:x+}%
				\(
					x_{\mathcal J^{k+1}}^{k+1}
				{}={}
					u_{\mathcal J^{k+1}}^k
				\)
				~and~
				\(
					x_{[N]\setminus\mathcal J^{k+1}}^{k+1}
				{}={}
					x_{[N]\setminus\mathcal J^{k+1}}^k
				\)
			\end{algorithmic}
		\end{algorithm}
		

		\subsection{Block-coordinate proximal point reformulation of \texorpdfstring{\cref{alg:Finito}}{Algorithm \ref*{alg:Finito}}}%
			\cref{alg:BCPP} presents a block coordinate (BC) proximal point algorithm  with the distance generating function \(\h*\).
			Note that in a departure from most of the existing literature on BC proximal methods that consider \emph{separable} nonsmooth terms (see \eg \cite{tseng2009coordinate,nesterov2012efficiency,beck2013convergence,bolte2014proximal,gao2020randomized}), here the nonsmooth function \(G\) in \eqref{eq:bigP} is nonseparable. It is shown in the next lemma that this conceptual algorithm is equivalent to the Bregman Finito/\allowbreak MISO \cref{alg:Finito}.
			
			\begin{lemma}[equivalence of \cref{alg:Finito,alg:BCPP}]\label{thm:equivBC}%
				As long as the same initialization parameters are chosen in the two algorithms, to any sequence \(\seq{\bm s^k,\s^k,z^k,\mathcal I^{k+1}}\) generated by \cref{alg:Finito} there corresponds a sequence \(\seq{\bm x^k,\bm u^k,\mathcal J^{k+1}}\) generated by \cref{alg:BCPP} (and viceversa) satisfying the following identities for all \(k\in\N\) and \(i\in [N]\):%
				\begin{enumerate}
				\item
					\(\mathcal I^{k+1}=\mathcal J^{k+1}\)
				\item\label{thm:equivBC:z}%
					\((z^k,\ldots,z^k)=\bm u^k\)
				\item\label{thm:equivBC:s}%
					\(s_i^k=\frac{1}{\gamma_i}\nabla \h_i(x_i^k)-\tfrac1N\nabla f_i(x_i^k)\)~ (or, equivalently, \(x_i^k=\nabla\conj{\h*_i}(s_i^k)\))
				\item\label{thm:equivBC:ts}%
					\(\tilde s^k=\sum_{i=1}^N\nabla\h*_i(x_i^k)\)
				\item\label{thm:equivBC:phi}%
					\(\varphi(z^k)=\Phi(\bm u^k)=\FBE(\bm x^k)-\D*(\bm u^k,\bm x^k)\)
				\item\label{thm:equivBC:L}%
					\(\FBE(\bm x^k)=\LL(z^k,\bm s^k)\)\bl, where \(\LL\) is as in \eqref{eq:L}.
				\end{enumerate}
				\begin{proof}
					Let the index sets \(\mathcal I^{k+1}\) and \(\mathcal J^{k+1}\) be chosen identically, \(k\in\N\).
					It follows from \cref{thm:Legendre:rangeT} that \(u_i^k=u_j^k\) for all \(k\in\N\) and \(i,j\in[N]\), with
					\begin{equation}\label{eq:vk}
						u_i^k
					{}={}
						\smash{
							\argmin_{w\in\R^n}\set{
								\textstyle
								g(w)
								{}+{}
								\sum_{i=1}^N{
									\tfrac{1}{\gamma_i}\h_i(w)
								}
								{}-{}
								\innprod*{
									\overbracket*{
										\textstyle
										\sum_{i=1}^N
											\nabla \h*_i(x_i^k)
									}^{\coloneqq v^k}
								}{w}
							}.
						}
					\end{equation}
					We now proceed by induction to show assertions \ref{thm:equivBC:z}, \ref{thm:equivBC:s}, and \ref{thm:equivBC:ts}.
					Note that the latter amounts to showing that \(v^k\) as defined in \eqref{eq:vk} coincides with \(\s^k\); by comparing \eqref{eq:vk} and the expression of \(z^k\) in \cref{state:Finito:z}, the claimed correspondence of \(\bm u^k\) and \(z^k\) as in assertion \ref{thm:equivBC:z} is then also obtained and, in turn, so is the identity in \ref{thm:equivBC:phi}.
					
					For \(k=0\) assertions \ref{thm:equivBC:s} and \ref{thm:equivBC:ts} hold because of the initialization of \(\s^0\) in \cref{alg:Finito} and of \(\bm x^0\) in \cref{alg:BCPP}; in turn, as motivated above, the base case for assertion \ref{thm:equivBC:z} also holds.
					Suppose now that the three assertions hold for some \(k\geq0\); then,
					\def\mycmd#1{\fillwidthof[c]{\nabla\h*_i(u_i^k)}{#1}}%
					\begin{align*}
					\textstyle
						v^{k+1}
					{}={}
						\sum_{i=1}^N\nabla\h*_i(x_i^{k+1})
					{}={} &
						\sum_{i\in\mathcal I^{k+1}}\mycmd{\nabla\h*_i(u_i^k)}
						{}+{}
						\;v^k\;
						{}-{}
						\sum_{i\in\mathcal I^{k+1}}\nabla\h*_i(x_i^k)
					\\
						\dueto{(induction)}
					{}={} &
						\sum_{i\in\mathcal I^{k+1}}{
							\mycmd{\nabla\h*_i(z^k)}
						}
						{}+{}
						\;\s^k\;
						{}-{}
						\sum_{i\in\mathcal I^{k+1}}{
							s_i^k
						}
					\\
					{}={} &
						\sum_{i\in\mathcal I^{k+1}}{
							\mycmd{s_i^{k+1}}
						}
						{}+{}
						\;\s^k\;
						{}-{}
						\sum_{i\in\mathcal I^{k+1}}{
							s_i^k
						}
					{}={}
						\s^{k+1},
					\end{align*}
					where the last two equalities are due to \cref{state:Finito:I,state:Finito:s+}.
					Therefore, \(v^{k+1}=\s^{k+1}\) and thus \(\bm u^{k+1}=(z^{k+1},\ldots,z^{k+1})\).
					It remains to show that
					\(
						s_i^{k+1}
					{}={}
						\tfrac{1}{\gamma_i}\nabla \h_i(x_i^{k+1})-\tfrac1N\nabla f_i(x_i^{k+1})
					\).
					For \(i\in\mathcal I^{k+1}\) this holds because of the update rule at \cref{state:Finito:s+} and the fact that \(x_i^{k+1}=u_i^k=z^k\) owing to \cref{state:BCPP:x+}.
					For \(i\notin\mathcal I^{k+1}\) this holds because \((x_i^{k+1},s_i^{k+1})=(x_i^k,s_i^k)\).
					Finally,
					\[
					\smash{
						\FBE(\bm x^k)
					{}\defeq{}
						\Phi(\bm u^k)
						{}+{}
						\D*(\bm u^k,\bm x^k)
					{}\overrel{\ref{thm:equivBC:phi}}{}
						\varphi(z^k)
						{}+{}
						\sum_{i=1}^N\D*_i(z^k,x_i^k)
					{}\overrel{\ref{thm:equivBC:s}}{}
						\varphi(z^k)
						{}+{}
						\sum_{i=1}^N\D*_i(z^k,\nabla\conj{\h*_i}(s_i^k)),
					}
					\vphantom{\sum}
					\]
					and the last term is \(\LL(z^k,\bm s^k)\) (cf. \cref{thm:D*,thm:inv1}), yielding assertion \ref{thm:equivBC:L}.
				\end{proof}
			\end{lemma}

	\section{Convergence analysis}\label{sec:conv}%
		The block coordinate interpretation of \cref{alg:Finito}  presented in \cref{sec:BC} plays a crucial role in the proposed methodology, and leads to a remarkably simple convergence analysis.
		In fact, many key facts can be established without confining the discussion to a particular sampling strategy.
		These preliminary results are presented in the next subsection and will be extensively referred to in the subsequent subsections that are instead devoted to a specific sampling strategy.

		\subsection{General sampling results}%
			Unlike classical analyses of BC proximal methods that employ the cost as a Lyapunov function (see \eg, \cite[\S11]{beck2017first}), here, the nonseparability of \(G\) precludes this possibility.
			To address this challenge, we instead employ the \emph{Bregman Moreau envelope} equipped with the distance generating function $\h*$ (see \eqref{eq:hatH}).
			Before showing its Lyapunov-type behavior for \cref{alg:BCPP}, we list some of its properties and its relation with the original problem.
			The proof is a simple consequence of \cref{thm:C0,thm:Moreau:basic} and the fact that \(\h*\) is a Legendre kernel with \(\dom\h*=\dom\h_1\times\dots\times\dom\h_N\) (cf. \cref{thm:Legendre}).
			
			\begin{lemma}[connections between \(\varphi+\indicator_{\overline C}\) and \(\FBE\)]%
				\bl{If \cref{ass:basic} holds, then:}%
				\begin{enumerate}
				\item
					\(\FBE\) is continuous on \(\dom\FBE=\interior\dom\h_1\times\dots\times\interior\dom\h_N\), in fact, locally Lipschitz if so is \(\nabla\h_i\) on \(\interior\dom\h_i\), \(i\in[N]\);
				\item \label{lem:env:prop:inf}
					\(
						\min_{\overline C}\varphi
					{}\leq{}
						\inf_C\varphi
					{}={}
						\inf\FBE
					\)
					and
					\(
						\argmin\FBE
					{}={}
						\set{(x^\star,\dots,x^\star)}[x^\star\in\argmin_C\varphi]
					\);
				\item \label{lem:envCoercive}
					\(\FBE\) is level bounded iff so is \(\varphi+\indicator_{\overline C}\).
				\end{enumerate}
			\end{lemma}
			
			\begin{lemma}[sure descent]\label{thm:sure}%
				Suppose that \cref{ass:basic} holds, and consider the iterates generated by \cref{alg:BCPP}.
				Then, \(\bm u^k=(u^k,\ldots,u^k)\) for some \(u^k\in C\) and \(\bm x^k\in C\times\dots\times C\subseteq\interior\dom\h*\) for every \(k\in\N\), and the algorithm is thus well defined.
				Moreover:
				\begin{enumerate}
				\item\label{thm:Igeq}%
					\(
						\FBE(\bm x^{k+1})
					{}\leq{}
						\FBE(\bm x^k)
						{}-{}
						\D*(\bm x^{k+1},\bm x^k)
					{}={}
						\FBE(\bm x^k){}
						{}-{}
						\sum_{i\in\mathcal J^{k+1}}\D*_i(u^k,x_i^k)
					\)
					for every \(k\in\N\);
					when \(\Phi\) is convex (\ie, when so is \(\varphi\)), then the inequality can be strengthened to
					\(
						\FBE(\bm x^{k+1})
					{}\leq{}
						\FBE(\bm x^k)
						{}-{}
						\D*(\bm x^{k+1},\bm x^k)
						{}-{}
						\D*(\bm u^k,\bm u^{k+1})
					\).
				\item\label{thm:decrease}%
					\(\seq{\FBE(\bm x^k)}\) monotonically decreases to a finite value \(\varphi_\star\geq\inf_C\varphi\geq\min_{\overline C}\varphi\).
				\item\label{thm:xdiff}%
					The sequence \(\seq{\D*(\bm x^{k+1},\bm x^k)}\) has finite sum (and in particular vanishes);
					the same holds also for \(\seq{\D*(\bm u^k, \bm u^{k+1})}\) when \(\Phi\) is convex (\ie, when so is \(\varphi\)).
				\item\label{thm:bounded}%
					If \(\varphi+\indicator_{\overline C}\) is level bounded, then \(\seq{\bm x^k}\) and \(\seq{\bm u^k}\) are bounded.
				\item\label{thm:x1}%
					If \(\dom\h_i\) is closed, a subsequence \(\seq{x_i^k}[k\in K]\) converges to a point \(x^\star\) iff so does \(\seq{x_i^{k+1}}[k\in K]\).%
				\item\label{thm:omega}%
					If \(C=\R^n\), then \(\FBE\) is \bl{constantly equal to \(\varphi_\star\) as above on the limit set of \(\seq{\bm x^k}\)}.%
				\end{enumerate}
				\begin{proof}
					It follows from \cref{thm:Legendre:rangeT} that \(u^k\in C\) holds for every \(k\in\N\).
					Notice that for every \(i\in[N]\) and \(k\in\N\), either \(x_i^k=x^{\rm init}\in C\) (by initialization), or there exists \(k_i\leq k\) such that \(x_i^k=z^{k_i}\in C\).
					It readily follows that \(\bm x^k\in C\times\dots\times C\subseteq\interior\dom\h=\interior\dom\h*\), hence that \(\prox_\Phi^{\h*}(\bm x^k)\neq\emptyset\) for all \(k\in\N\) by \cref{thm:Legendre:rangeT}, whence the well definedness of the algorithm.
					We now show the numbered claims.
					\begin{proofitemize}
					\item\ref{thm:Igeq}~
						It follows from \cref{thm:leq,thm:leq_convex} that
						\(
							\FBE(\bm x^{k+1})
						{}\leq{}
							\Phi(\bm u^k)+\D*(\bm u^k,\bm x^{k+1})
							{}-{}
							c_k
						\),
						where \(c_k\geq0\) can be taken as \(c_k=\D*(\bm u^k,\bm u^{k+1})\) when \(\Phi\) is convex.
						Therefore,
						\begin{align*}
							\FBE(\bm x^{k+1})
						{}\leq{} &
							\Phi(\bm u^k)+\D*(\bm u^k,\bm x^{k+1})
							{}-{}
							c_k
						{}={}
							\FBE(\bm x^k)
							{}-{}
							\D*(\bm u^k,\bm x^k)
							{}+{}
							\D*(\bm u^{k},\bm x^{k+1})
							{}-{}
							c_k
						\\
						{}\overrel*{\ref{thm:3P}}{} &
							\FBE(\bm x^k)
							{}-{}
							\D*(\bm x^{k+1},\bm x^k)
							{}-{}
							\innprod*{\bm u^k-\bm x^{k+1}}{
								\nabla\h*(\bm x^{k+1})-\nabla\h*(\bm x^k)
							}
							{}-{}
							c_k.
						\end{align*}
						The claim follows by noting that the inner product is zero:
						\begin{align*}
							\innprod*{\bm u^k-\bm x^{k+1}}{
								\nabla\h*(\bm x^{k+1})-\nabla\h*(\bm x^k)
							}
						{}={}
							\sum*_{j\in[N]}{
								\innprod*{
									u^k
									{}-{}
									\underbracket*{x_j^{k+1}}_{
										=u^k
										\mathrlap{\text{ for }j\in\mathcal J^{k+1}}
									}
								}{
									\nabla\h*_j(
										\overbracket*{x_j^{k+1}}^{
											=x_j^k
											\mathrlap{\text{ for }j\notin\mathcal J^{k+1}}
										}
									)
									{}-{}
									\nabla\h*_j(x_j^k)
								}
							}=0.
						\end{align*}
					\item\ref{thm:decrease}~
						Monotonic decrease of \(\seq{\FBE(\bm x^k)}\) follows from assertion \ref{thm:Igeq}.
						This ensures that the sequence converges to some value \(\varphi_\star\), bounded below by \(\min_{\overline C}\varphi\) in light of \cref{lem:env:prop:inf}.
					\item\ref{thm:xdiff}~
						It follows from assertion \ref{thm:Igeq} that
						\[\textstyle
							\sum_{k\in\N}\D*(\bm x^{k+1},\bm x^k)
						{}\leq{}
							\FBE(\bm x^0)-\inf\FBE
						{}\leq{}
							\FBE(\bm x^0)-\inf_{\overline C}\varphi
						{}<{}
							\infty
						\]
						owing to \cref{lem:env:prop:inf} and \cref{ass:argmin}.
						When \(\varphi\) is convex, the tighter bound in assertion \ref{thm:Igeq} yields the similar claim for \(\seq{\D*(\bm u^k,\bm u^{k+1})}\).
					\item\ref{thm:bounded}~
						It follows from assertion \ref{thm:decrease} that the entire sequence \(\seq{\bm x^k}\) is contained in the sublevel set \(\set*{\bm w}[\FBE(\bm w)\leq\FBE(\bm x^0)]\), which is bounded provided that \(\varphi+\indicator_{\overline C}\) is level bounded as shown in \cref{lem:envCoercive}.
						In turn, boundedness of \(\seq{\bm u^k}\) then follows from local boundedness of \(\T=\prox_\Phi^{\h*}\), cf. \eqref{eq:Tprox} and \cref{thm:osc}.
					\item\ref{thm:x1}~
						Follows from \cref{thm:Dboundary}, since \(x_i^k\in\interior\dom\h_i=\interior\dom\h*_i\) for every \(k\) (with equality owing to \cref{thm:Legendre}), and \(\D*_i(x_i^{k+1},x_i^k)\to0\) by assertion \ref{thm:xdiff}.
					\item\ref{thm:omega}~
						Follows from assertion \ref{thm:decrease} and the continuity of \(\FBE\), see \cref{thm:C0}.
					\qedhere
					\end{proofitemize}
				\end{proof}
			\end{lemma}
			
			In conclusion of this subsection we provide an overview of the ingredients that are needed to show that the limit points of the sequence \(\seq{z^k}\) generated by \cref{alg:Finito} are stationary for problem \eqref{eq:P}.
			As will be shown in \cref{thm:subseq}, these amount to the vanishing of the \emph{residual} \(\D*(\bm u^k,\bm x^k)\) together with some assumptions on the distance-generating functions \(\h_i\).
			For the iterates of \cref{alg:Finito}, this translates to
			\(
				\BregmanD_{\conj{\h*_i}}(s_i^k,\nabla\h*_i(z^k))
			{}\to{}
				0
			\)
			for all indices \(i\in[N]\), indicating that all vectors \(s_i^{k+1}\) in the table should be good estimates of
			\(
				\nabla\h*_i(z^{k+1})
			{}={}
				\frac{1}{\gamma_i}\nabla\h_i(z^{k+1})
				{}-{}
				\frac1N\nabla f_i(z^{k+1})
			\),
			as opposed to
			\(
				\frac{1}{\gamma_i}\nabla\h_i(z^k)
				{}-{}
				\frac1N\nabla f_i(z^k)
			\)
			and for the indices in \(\mathcal I^{k+1}\) only (cf. \cref{state:Finito:s+}).
			As a result, we may view this property as jointly having \(z^k-z^{k+1}\) vanish, desirable if any convergence of \(\seq{z^k}\) is expected, and the fact that a consensus is eventually reached among the sampled blocks.
			
			In line with any result in the literature we are aware of, a complete convergence analysis for nonconvex problems will ultimately require \(C=\R^n\).
			For convex problems, that is, when the cost function \(\varphi\) is convex without any among \(f_i\) and \(g\) being necessarily so, the following requirement will instead suffice to our purposes in the randomized sampling setting of \eqref{eq:RS}.
			
			\begin{assumption}[requirements on the distance-generating functions]\label{ass:d}%
				For \(i\in[N]\), \(\dom\h_i\) is closed, and whenever \(\interior\dom\h_i\ni z^k\to z\in\boundary\dom\h_i\) it holds that \(\D_i(z,z^k)\to0\).
			\end{assumption}
			
			\begin{remark}
				\begin{enumerate}
				\item\label{thm:fulldom}%
					\Cref{ass:d} is vacuously satisfied when \(\dom\h_i=\R^n\), having \(\boundary\R^n=\emptyset\).
			\bl{%
				\item\label{rem:AssII:ex}%
					While \Cref{ass:d} always holds on \(\R\), it may fail in higher dimensions \cite[Ex. 7.32]{bauschke1997legendre}.%
			}%
				\item\label{thm:h*}%
					For any \(i\in[N]\), function \(\h_i\) complies with \cref{ass:d} iff so does \(\h*_i\), owing to the inequalities
					\(
						\tfrac{N-\gamma_i L_{f_i}}{N\gamma_i}\D_i
					{}\leq{}
						\D*_i
					{}\leq{}
						\tfrac{N+\gamma_i L_{f_i}}{N\gamma_i}\D_i
					\)
					(cf. \eqref{eq:hhat_LS}).
				\hfill\qedsymbol
				\end{enumerate}\let\qedsymbol\relax
			\end{remark}
			
			\begin{lemma}[subsequential convergence recipe]\label{thm:subseq}%
				Suppose that \cref{ass:basic} holds, and consider the iterates generated by \cref{alg:Finito}.
				Let \(x_i^k=\nabla\conj{\h*_i}(s_i^k)\) and \(z^k=u^k\) be the corresponding iterates generated by \cref{alg:BCPP} as in \cref{thm:equivBC}, and suppose that
				\begin{enumeratass}
				\item\label{ass:res}%
					\(\D*(\bm u^k,\bm x^k)\to0\) (or equivalently,
					\(
						\BregmanD_{\conj{\h*_i}}(s_i^k,\nabla\h*_i(z^k))
					{}\to{}
						0
					\),
					\(i\in[N]\)).
				\end{enumeratass}
				Then, letting \(\varphi_\star\) be as in \cref{thm:decrease}, the following hold:
				\begin{enumerate}
				\item\label{thm:cost}%
					\(\varphi(z^k)=\Phi(\bm u^k)\to\varphi_\star\) as \(k\to\infty\).
				\item\label{thm:xzlim}%
					If \(\dom\h_i\) is closed, \(i\in[N]\), then having
					(a) \( \seq{z^k}[k\in K]\to z\), ~
					(b) \(\seq{x_i^k}[k\in K]\to z\) \(\exists i\in[N]\), ~
					and
					(c) \(\seq{z^{k+1},x_i^{k+1}}[k\in K]\to(z,z)\) \(\forall i\in[N]\),
					are all equivalent conditions.
					In particular, if \(\seq{z^k}\) is bounded (\eg when \(\varphi+\indicator_{\overline C}\) is level bounded), then \(\|z^{k+1}-z^k\|\to0\) holds, and the set of its limit points, be it \(\omega\), is thus nonempty, compact, and connected.%
				\item\label{thm:constant}%
					Under \cref{ass:d}, \(\varphi\equiv\varphi_\star\) on \(\omega\) (the set of limit points of \(\seq{z^k}\)).
				\item\label{thm:stationary}%
					If \(C=\R^n\), then every \(z^\star\in\omega\) is stationary for \eqref{eq:P}.
				\end{enumerate}
				\begin{proof}
				\cref{ass:res} can be written as
					\(
						\D*_i(z^k,\nabla\conj{\h*_i}(s_i^k))
					{}\to{}
						0
					\),
					\(i\in[N]\).
					In turn, by the conjugate identity \bl{in} \cref{thm:D*}, the equivalent expression in terms of $s_i^k$ and $z^k$ is obtained.%
					\begin{proofitemize}
					\item\ref{thm:cost}~
						As shown in \cref{thm:decrease}, \(\seq{\FBE(\bm x^k)}\) monotonically decreases to \(\varphi_\star\).
						In turn, \cref{thm:equivBC:phi} and \cref{ass:res} then imply that \(\varphi(z^k)=\Phi(\bm u^k)\) converges to \(\varphi_\star\).
					\item\ref{thm:xzlim}~
						The equivalences owe to \cref{thm:Dboundary,thm:x1} (as $\dom\h*_i=\dom\h_i$), and imply \(\|z^{k+1}-z^k\|\to0\) if \(\seq{z^k}\) is bounded.
						The claim on \(\omega\) then follows from \cite[Rem. 5]{bolte2014proximal}.%
					\item\ref{thm:constant}~
						Let \(z^\star\in\omega\) be fixed, and let \(\seq{z^k}[k\in K]\) be a subsequence converging to \(z^\star\).
						Assertion \ref{thm:xzlim} ensures that \(\seq{\bm x^k}[k\in K]\to\bm z^\star\coloneqq(z^\star,\dots,z^\star)\), hence
						\[
							\varphi_\star
						{}\xleftarrow[k\in K]{\text{\ref{thm:decrease}}}{}
							\FBE(\bm x^k) {}-{} \D*(\bm z^\star,\bm x^k)
						{}\overrel[\leq]{\ref{thm:leq}}{}
							\Phi(\bm z^\star)
						{}\overrel[\leq]{lsc}{}
							\liminf_{k\in K}\Phi(\bm u^k)
						{}\overrel{\ref{thm:cost}}{}
							\varphi_\star,
						\]
						where \cref{ass:d} is used in the first limit.
					\item\ref{thm:stationary}~
						Suppose that \(C=\R^n\) and \(\seq{z^k}[k\in K]\to z^\star\) for some infinite \(K\subseteq\N\) and \(z^\star\in\R^n\), so that, by virtue of assertion \ref{thm:xzlim}, \(\seq{\bm x^k,\bm x^{k+1}}[k\in K]\to(\bm z^\star,\bm z^\star)\).
						Since \((z^k,\ldots,z^k)=\bm u^k\in\prox_\Phi^{\h*}(\bm x^k)\), the osc of \(\prox_\Phi^{\h*}\) (\cref{thm:osc}) ensures that
						\(
							\bm z^\star
						{}\in{}
							\prox_{\Phi}^{\h*}(\bm z^\star)
						\),
						hence \(0\in\hat\partial\Phi(\bm z^\star)\) owing to \cref{thm:OC}.
						By invoking \cref{thm:equivP:stationary} we conclude that \(z^\star\) is stationary for \eqref{eq:P}.
					\qedhere
					\end{proofitemize}
				\end{proof}
			\end{lemma}
			

		\subsection{Randomized sampling rule\texorpdfstring{ \eqref{eq:RS}}{}}\label{sec:RS}%

			The analysis for the randomized case dealt in this section will make use of the following result, known as the Robbins-Siegmund supermartingale theorem, and stated here in simplified form following \cite[Prop. 2]{bertsekas2011incremental}.
			
			\begin{fact}[supermartingale convergence theorem {\cite{robbins1985convergence}}]\label{thm:Robbins}%
				For \(k\in\N\), let \(\xi_k\) and \(\eta_k\) be random variables, and \(\mathcal F_k\subseteq\mathcal F_{k+1}\) be sets of random variables such that
				\begin{enumeratass}
				\item
					\(\mathcal F_k\subseteq\mathcal F_{k+1}\);
				\item
					\(0\leq\xi_k,\eta_k\) are functions of the random variables in \(\mathcal F_k\);
				\item
					\(
						\E[]{~\xi_{k+1}~}[~\mathcal F_k~]
					{}\leq{}
						\xi_k-\eta_k
					\).
				\end{enumeratass}
				Then, almost surely, \(\sum_{k\in\N}\eta_k<\infty\) and \(\xi_k\) converges to a (nonnegative) random variable.
			\end{fact}
			
			The sets \(\mathcal F_k\) in the above formulation will represent the information available at iteration \(k\), and the notation \(\E{{}\cdot{}}\) will be used as a shorthand for \(\E[]{{}\cdot{}}[~\mathcal F_k~]\).
			
			\begin{theorem}[subsequential convergence of \cref{alg:Finito} with randomized \bl{rule} \eqref{eq:RS}]\label{thm:RS:subseq}%
				Suppose that \cref{ass:basic} holds.
				Then, denoting \(p_{\rm min}=\min_ip_i\), the iterates generated by \cref{alg:Finito} with indices selected according to the randomized rule \eqref{eq:RS} satisfy
				\begin{equation}\label{eq:RS:descent}
					\E{\LL(z^{k+1},\bm s^{k+1})}
				{}\leq{}
					\LL(z^k,\bm s^k)
					{}-{}
					p_{\rm min}\sum_{i=1}^N\D**_i(s_i^k,\nabla\h*_i(z^k))
				\quad
					\forall k\in\N,
				\end{equation}
				where \(\LL\) is as in \eqref{eq:L} (and satisfies \(\LL(z^k,\bm s^k)=\FBE(\bm x^k)\), cf. \cref{thm:equivBC:L}).
				Moreover, letting \(\omega\) denote the set of limit points of \(\seq{z^k}\)\bl, the following assertions hold almost surely:
				\begin{enumerate}
				\item\label{thm:RS:res}%
					The sequence \(\seq{\D**_i(s_i^k,\nabla\h*_i(z^k))}\) has finite sum (and in particular vanishes), $i\in[N]$.
				\item\label{thm:RS:cost}%
					The sequence $\seq{\varphi(z^k)}$ converges to the finite value $\varphi_\star\leq\varphi(x^{\rm init})$ of \cref{thm:decrease}.
				\item\label{thm:RS:const}%
					If \cref{ass:d} is satisfied, then \(\varphi\equiv\varphi_\star\) on \(\omega\).
				\item\label{thm:RS:stationary}%
					If \(C=\R^n\), then \(0\in\hat\partial\varphi(z^\star)\) for every \(z^\star\in\omega\).
				\end{enumerate}
				When \(\varphi\) is convex (without \(g\) or any \(f_i\) necessarily being so) and \(\dom h_i\) is closed, \(i\in[N]\),
				\begin{enumerate}[resume]
				\item\label{thm:RS:convex:cost}%
					\(\seq{\varphi(z^k)}\) converges to \(\min_{\overline C}\varphi\) with \(\bl{\min_{\ell=0,\ldots k}\set*{\varphi(z^\ell)}}-\min_{\overline C}\varphi \leq o(\nicefrac1k)\);
				\item\label{thm:RS:convex:stationary}%
					the limit points of \(\seq{z^k}\) all belong to \(\argmin_{\overline C}\varphi\);
				\item\label{thm:RS:convex:seq}%
					if either \cref{ass:d} holds and \(\varphi+\indicator_{\overline C}\) is level bounded, or \(C=\R^n\),
					then \(\seq{z^k}\) and \(\seq{\nabla \h*_i^*(s_i^k)}\), \(i\in[N]\), converge to the same point in \(\argmin_{\overline C}\varphi\).%
				\end{enumerate}
				\begin{proof}
					By \cref{thm:equivBC}, we \bl{will} consider the simpler setting of \cref{alg:BCPP}.
					We have%
					\begin{align*}
						\E{\LL(z^{k+1},\bm s^{k+1})}
					{}\overrel{\ref{thm:equivBC:L}}{} &
						\E{\FBE(\bm x^{k+1})}
					{}\overrel[\leq]{\ref{thm:Igeq}}{}
						\E{\FBE(\bm x^k)
						{}-{}
						\smash{\sum_{i\in\mathcal I^{k+1}}}\D*_i(u^k,x_i^k)}
					\\
					{}={} &
						\FBE(\bm x^k)
						{}-{}
						\smash{\sum_{i=1}^N}p_i\D*_i(u^k,x_i^k)
						\vphantom{\sum^N}
					\\
					{}\leq{} &
						\FBE(\bm x^k)
						{}-{}
						p_{\rm min}\D*(\bm u^k,\bm x^k)
					{}\overrel{\ref{thm:equivBC:L}, \ref{thm:D*}}{}
						\LL(z^k,\bm s^k)
						{}-{}
						p_{\rm min}
						\smash{\sum_{i=1}^N}\D**_i(s_i^k,\nabla\h*_i(z^k))\bl,
					\end{align*}
					which is \eqref{eq:RS:descent}.
					We thus infer from \cref{thm:Robbins} that \(\seq{\D*(\bm u^k,\bm x^k)}\) has almost surely finite sum, and the proof of assertions \ref{thm:RS:res}--\ref{thm:RS:stationary} then follows from \cref{thm:subseq}.
			
					\bl{%
						In what follows, suppose that \(\varphi\) is convex and \(\dom h_i\) is closed, \(i\in[N]\), so that
						\(
							\overline C
						{}={}
							\bigcap_{i=1}^N\overline{\interior\dom\h_i}
						{}={}
							\bigcap_{i=1}^N\dom\h_i
						\)
						\cite[Prop. 1.3.8]{bertsekas2015convex}, and in particular \(\D*_i(y,x)<\infty\) holds for any \((y,x)\in\dom\h_i\times\interior\dom\h_i\supseteq\overline C\times C\).
					}%
					\begin{proofitemize}
					\item\ref{thm:RS:convex:cost}~
						For any \(x^\star\in\argmin_{\overline C}\varphi\), so that \(\bm x^\star\coloneqq(x^\star,\ldots,x^\star)\in\argmin_{\overline C\times\cdots\times\overline C}\Phi\) (cf. \cref{thm:equivP:argmin}), the three-point identity (\cref{thm:3P}), convexity of \(\Phi\) (\cref{thm:equivP:convex}) and the inclusion
						\(
							\nabla\h*(\bm x^k)-\nabla\h*(\bm u^k)
						{}\in{}
							\hat\partial\Phi(\bm u^k)
						\)
						(\cref{thm:OC}) give the bound
						\begin{align*}
							\D*(\bm x^\star,\bm u^k)
						{}={} &
							\D*(\bm x^\star,\bm x^k)
							{}-{}
							\D*(\bm u^k,\bm x^k)
							{}+{}
							\innprod*{\nabla \h*(\bm x^k)-\nabla \h*(\bm u^k)}{\bm x^\star-\bm u^k}
						\\
						\numberthis\label{eq:Phi_convex}
						{}\leq{} &
							\D*(\bm x^\star,\bm x^k)
							{}-{}
							\D*(\bm u^k,\bm x^k)
							{}+{}
							\Phi(\bm x^\star)
							{}-{}
							\Phi(\bm u^k).
						\end{align*}
						Next,
						\begin{align*}
							\E{\sum_{i=1}^Np_i^{-1}\D*_i(x^\star,x_i^{k+1})}
						{}={} &
							\sum_{i=1}^Np_i^{-1}\left(
								\overbracket*{\P{\mathcal I^{k+1}\ni i}}^{p_i}
								\D*_i(x^\star,u^k)
								{}+{}
								\overbracket*{\P{\mathcal I^{k+1}\not\ni i}}^{1-p_i}
								\D*_i(x^\star,x^k)
							\right)
						\\
						\numberthis\label{eq:RS:convex:prob}
						{}={} &
							\D*(\bm x^\star,\bm u^k)
							{}+{}
							\smash{\sum_{i=1}^N(1-p_i)p_i^{-1}\D*_i(x^\star,x_i^k)}
						\\
							\dueto{\eqref{eq:Phi_convex}}
						{}\leq{} &
							\D*(\bm x^\star,\bm x^k)
							{}-{}
							\D*(\bm u^k,\bm x^k)
							{}+{}
							\Phi(\bm x^\star)
							{}-{}
							\Phi(\bm u^k)
							{}+{}
							\sum_{i=1}^N(p_i^{-1}-1)\D*_i(x^\star,x_i^k)
						\\
							\dueto{\cref{thm:equivP,thm:equivBC:phi}}
						{}={} &
							\smash{\sum_{i=1}^Np_i^{-1}\D*_i(x^\star,x_i^k)}
							{}-{}
							\D*(\bm u^k,\bm x^k)
							{}-{}
							\bigl(
								\varphi(\bm z^k)
								{}-{}
								\min_{\overline C}\varphi
							\bigr),
						\end{align*}
						where \(\bm u^k=(u^k,\ldots,u^k)\).
						From \cref{thm:Robbins} we conclude that
						\begin{equation}\label{eq:RS:convex:RobbinsConseq}
							\smash{\sum_{k\in\N}}\D*(\bm u^k,\bm x^k)<\infty
						\quad\text{and}\quad
							\smash{\sum_{k\in\N}}\bigl(\varphi(z^k)-\min_{\overline C}\varphi\bigr) \eqqcolon c
						{}<{}
							\infty
						\end{equation}
						almost surely, and
						\begin{equation}\label{eq:RS:convex:dconvas}
							\smash{\sum_{i=1}^N}p_i^{-1}\D*_i(x^\star,x_i^k)
						\quad
							\text{converges a.s. for any \(x^\star\in\argmin_{\overline C}\varphi\).}
						\end{equation}
						It now follows from \eqref{eq:RS:convex:RobbinsConseq} that \(\varphi(z^k)\) converges a.s. to \(\min_{\overline C}\varphi\).
						\bl{Moreover, since the sequence \(\seq{\min_{\ell=0,\ldots k}\varphi(z^\ell)}\) is nonincreasing, \eqref{eq:RS:convex:RobbinsConseq} also yields the claimed rate}.
					\item\ref{thm:RS:convex:stationary}~
						Suppose that \(\seq{z^k}[k\in K]\to z^\star\).
						Then, \(\seq{\bm u^k}[k\in K]\to\bm u^\star\) for \(\bm u^k=(z^k,\ldots,z^k)\) and \(\bm u^\star=(z^\star,\ldots,z^\star)\).
						Notice that \(z^\star\in\overline C\), since \(z^k\in C\) for all \(k\) (cf. \cref{thm:sure}).
						We have
						\begin{equation}
							\min_{\overline C}\varphi
						{}\overrel{\ref{thm:equivP:argmin}}{}
							\min_{\overline C\times\cdots\times\overline C}\Phi
						{}\leq{}
							\Phi(\bm u^\star)
						{}\overrel[\leq]{lsc}{}
							\liminf_{k\in K}\Phi(\bm u^k)
						{}\overrel{\ref{thm:equivBC:phi}}{}
							\liminf_{k\in K}\varphi(z^k)
						{}\overrel{\ref{thm:RS:convex:cost}}{}
							\min_{\overline C}\varphi.
						\end{equation}
						Therefore, \(\bm u^\star\) is a minimizer of \(\Phi\) on \(\overline C\times\cdots\times\overline C\), and thus \(z^\star\) is a minimizer of \(\varphi\) on \(\overline C\) by virtue of \cref{thm:equivP:argmin}.
					\item\ref{thm:RS:convex:seq}~
						If \(\varphi+\indicator_{\overline C}\) is level bounded, then by \cref{thm:bounded} \(\seq{\bm x^k}\) and \(\seq{\bm u^k}\) are bounded.
						Alternatively, if \(C=\R^n\), then boundedness of the former sequence follows from \cref{thm:Dhcoercive}, \eqref{eq:RS:convex:dconvas} \bl{and \cref{ass:argmin}}, and in turn that of the latter from \eqref{eq:RS:convex:RobbinsConseq}.
						In either cases \cref{ass:d} holds, as discussed in \cref{thm:fulldom}.
						Boundedness of the sequences ensures the existence of \(K\subseteq\N\), \(z^\star\) and \(\bm u^\star\) as in the proof of assertion \ref{thm:RS:convex:stationary}.
						The vanishing of \(\D*(\bm u^k,\bm x^k)\) shown in \eqref{eq:RS:convex:RobbinsConseq} implies through \cref{thm:xzlim} that \(\seq{\bm x^k}\) and \(\seq{\bm u^k}\) have same limit points, and that \(\seq{\bm x^k}[k\in K]\to\bm u^\star\).
						In turn, \(\seq{\sum_{i=1}^Np_i^{-1}\D*_i(u^\star,x_i^k)}[k\in K]\to0\) \bl{holds by} \cref{ass:d}.
						Hence, since the entire sequence is convergent (by \eqref{eq:RS:convex:dconvas}) we have \(\seq{\sum_{i=1}^Np_i^{-1}\D*_i(u^\star,x_i^k)}\to0\), which by \cref{thm:Dboundary} implies \(\seq{x_i^k}\to u^\star\), \(i\in[N]\).
						As discussed above, this implies that \(\seq{\bm u^k}\to\bm u^\star\), and the identity \(\bm u^k=(z^k,\ldots,z^k)\) of \cref{thm:equivBC:z} yields the claimed convergence.%
					\qedhere
					\end{proofitemize}
				\end{proof}
			\end{theorem}
			\bl{%
			In \cref{thm:RS:convex:seq} the assumption that \(\varphi+\indicator_{\overline C}\) is level bounded can be relaxed by instead requiring that for every $v\in\dom \h_i$ and $\alpha\in\R$, the level set $\set{w\in\interior\dom\h_i}[\D_i(v,w)\leq\alpha]$ is bounded, as this would suffice to ensure boundedness of the sequences.
			In fact, together with the closed-domain requirement this is a standing assumption in many works dealing with Bregman distances, specifically those involving \emph{Bregman functions ``with zone \(S\)''} (\(S\) being the interior of the domain), see \eg, \cite{solodov2000inexact}.
			}%

			We conclude this subsection with an analysis of the strongly convex case, in which linear convergence (in expectation) will be shown.
			Remarkably, strong convexity of the cost function \(\varphi\) alone will suffice, without imposing any such requirement on the individual terms \(f_i\) or \(g\) which, in fact, are even allowed to be nonconvex.
			
			\begin{theorem}[linear convergence with randomized \bl{rule} \eqref{eq:RS} for strongly convex problems]\label{thm:RS:strconvex}%
				Consider the iterates of \cref{alg:Finito}.
				Additionally to \cref{ass:basic}, suppose that
				\begin{enumeratass}
				\item\label{ass:RS:strconvex:varphi}%
					\(\varphi\) is \(\mu_\varphi\)-strongly convex;
				\item\label{ass:RS:strconvex:h}%
					\(\h_i\) has locally Lipschitz gradient on the whole space \(\R^n\), \(i\in[N]\), (hence \(C=\R^n\)).
				\end{enumeratass}
				Let \(\U\) be a convex compact set containing \(x^{\rm init}\) and the sequence \(\seq{z^k}\), and let \(\ell_{\h_i,\U}\) be a Lipschitz modulus for \(\nabla\h_i\) on \(\U\), \(i\in[N]\).\footnote{%
					\(\U\) exists by \cref{thm:bounded}, owing to strong convexity and consequent level boundedness of \(\varphi\).
					For \(i\in[N]\), a finite \(\ell_{\h_i,\U}\) then exists because of \cref{ass:RS:strconvex:h} and since \(\U\subset C\) (as opposed to \(\U\subseteq\overline C\)), having \(C=\R^n\).%
				}
				Let \(x^\star=\argmin\varphi\), \(\varphi_\star=\min\varphi\), and%
				\begin{equation}\label{eq:RS:strconvex:c}
					c_{\U}
				{}={}
					\frac{\min_ip_i}{
						1
						{}+{}
						\frac{1}{\mu_\varphi}
						\sum_i{
							\bigl(
								\frac{\ell_{\h_i,\U}}{\gamma_i}
								{}-{}
								\frac{\sigma_{f_i,\U}}{N}
							\bigr)
						}
					},
				\end{equation}
				where \(\sigma_{f_i,\U}\geq -L_{f_i}\ell_{\h_i,\U}\) is a (weak) convexity modulus of \(f_i\) on \(\U\).\footnote{%
					\(f_i\) is \(\sigma_{f_i,\U}\)-weakly convex on \(\U\) if \(f_i-\frac{\sigma_{f_i,\U}}{2}\|{}\cdot{}\|^2\) is convex on \(\U\), thus coinciding with convexity (resp. \(\sigma_{f_i,\U}\)-strong convexity) on \(\U\) when \(\sigma_{f_i,\U}\geq0\) (resp. \(\sigma_{f_i,\U}>0\)).
					The lower bound on \(\sigma_{f_i,\U}\) owes to \cref{thm:fC11}.%
				}
				Then, for all \(k\in\N\)
				\begin{align}\label{eq:RS:strconvex:Qlinear}
					\E{\LL(z^{k+1},\bm s^{k+1})-\varphi_\star}
				{}\leq{} &
					(1-c_{\U})
					\bigl(\LL(z^k,\bm s^k)-\varphi_\star\bigr),
				\quad\text{and}
				\\
				\label{eq:RS:strconvex:Rlinear}
					\E[]{\tfrac{\mu_\varphi}{2}\|z^k-x^\star\|^2}
				{}\leq{}
					\E[]{\varphi(z^k)-\varphi_\star}
				{}\leq{} &
					(1-c_{\U})^k
					\bigl(\varphi(x^{\rm init})-\varphi_\star\bigr).
				\end{align}
				\begin{proof}
					\bl{See \cref{proof:thm:RS:strconvex}.}
				\end{proof}
			\end{theorem}
			
			In the Euclidean case, \(h_i=\tfrac12\|{}\cdot{}\|^2\) \bl{has \(1\)-Lipschitz gradient} on \(\R^n\), and the results of \cref{thm:RS:strconvex} hold with \(\U=\R^n\) and \(\ell_{\h_i}=1\), \bl{improving} those of \cite[Cor. 3.3]{latafat2021block} (limited to the Euclidean case) both by providing tighter rates and by relaxing (strong) convexity assumptions on individual \(f_i\).
			For the uniform sampling strategy \(p_i=\nicefrac1N\) the rate \(1-O(1/N)\) is obtained.
			The same arguments still hold for the Bregman extension of \cref{alg:Finito} dealt in this paper, as long as each \(\h_i\) is Lipschitz differentiable.
			This fact is stated in the following corollary, where, by using the fact that \(\mu_\varphi\geq\tfrac1N\sum_{i=1}^N\sigma_{f_i}\) under a convexity assumption on \(g\), a simplified expression for the constant \(c\) in \eqref{eq:RS:strconvex:c} is obtained.
			We remark that in the Euclidean case a variant of SVRG \cite{johnson2013accelerating} has also been studied in \cite{allenzhu2016improved} under similar assumptions.
			
			\begin{corollary}[global linear rate]\label{thm:RS:strconvex:O(1/N)}%
				Additionally to \cref{ass:basic}, suppose that
				\begin{enumeratass}
				\item
					\(g\) is convex, and \(f\coloneqq\frac1N\sum_if_i\) is \(\mu_f\)-strongly convex (yet each \(f_i\) can be nonconvex);%
				\item
					\(\nabla\h_i\) is Lipschitz on \(\R^n\) (hence so is \(f_i\) with modulus \(\ell_{f_i}\), cf. \cref{thm:fC11}), \(i\in[N]\).%
				\end{enumeratass}
				\bl{%
					Set \(\gamma_i=\nicefrac{\alpha N}{L_{f_i}}\) with \(\alpha\in(0,1)\) and
					\(
						\kappa_f
					{}\coloneqq{}
						\frac{\frac1N\sum_i\ell_{f_i}}{\mu_f}
					\).
					Then, \eqref{eq:RS:strconvex:Qlinear} and \eqref{eq:RS:strconvex:Rlinear} hold with
					\(
						c_{_{\R^n}}
					{}\geq{}
						\frac{\alpha\min_ip_i}{\kappa_f}
					\).
				}%
			\end{corollary}

		\subsection{Essentially cyclic sampling rule\texorpdfstring{ \eqref{eq:ECS}}{}}\label{Sec:ECS}%
			The convergence results in this subsection require convexity of the nonsmooth term \(g\) and local strong convexity and smoothness of \(\h_i\) (as is the case when \(h_i\in\C^2(\R^n)\) with \(\nabla^2\h_i\succ0\)).
			The proof of subsequential convergence is an adaptation of that of \cite[Thm. 2.8]{latafat2021block}. 
			
			\begin{theorem}[subsequential convergence with essentially cyclic \bl{rule} \eqref{eq:ECS}]\label{thm:ECS:subseq}%
				Additionally to \cref{ass:basic}, assume that \(g\) is convex, \(C=\R^n\), and that either \bl{one of the following assumptions holds:%
				\begin{enumerator}
				\item\label{ass:ECS:subseq:global}%
				{\color{black}%
					(either) each \(h_i\) is strongly convex and Lipschitz differentiable,%
				}%
				\item\label{ass:ECS:subseq:local}%
				{\color{black}%
					(or) \(\varphi\) is level bounded and each \(\h_i\) is locally strongly convex and locally Lipschitz differentiable.%
				}%
				\end{enumerator}}%
				Then, all the claims in \cref{thm:RS:res,thm:RS:cost,thm:RS:const,thm:RS:stationary} hold surely.
				\begin{proof}
					\bl{See \cref{proof:thm:ECS:subseq}.}
				\end{proof}
			\end{theorem}
			\bl{%
				We remark that linear convergence in the strongly convex case may be obtained in a similar fashion to \cite[Thm. 2.9]{latafat2021block} and \cite[Thm. 3.9]{beck2013convergence}.
				This however results in a rate more conservative than the one obtained for the randomized case, cf. \cite[Eq. (2.21)]{latafat2021block}, which is not consistent with what is observed in practice.
				Although a refined analysis not relying on conservative triangle inequalities may be possible, this direction is not investigated here.
			}%
			
			Our next goal is to establish global (and linear) convergence results without convexity assumptions on \(f_i\) or their sum.
			To this end, we leverage the  Kurdyka-\L ojasiewicz property \cite{lojasiewicz1993geometrie,kurdyka1998gradients}, which has become the standard tool in the analysis of nonconvex proximal methods, and most notably holds for the class of semialgebraic functions \cite{bolte2007clarke,bolte2007lojasiewicz,attouch2009convergence,attouch2010proximal,attouch2013convergence,bolte2014proximal}.
			
			\begin{definition}[KL property with exponent \(\theta\)]\label{def:KL}%
				A proper lsc function \(\func{q}{\R^n}{\Rinf}\) has the \DEF{Kurdyka-{\L}ojasiewicz} (KL) property with exponent \(\theta\in(0,1)\) if for every \(\bar w\in\dom\partial q\) there exist \(\varepsilon,\eta,\varrho>0\) such that
				\(
					\psi'(q(w)-q(\bar w))\dist(0,\partial q(w))\geq 1
				\)
				holds for every \(w\) satisfying \(\|w-\bar w\|<\varepsilon\) and \(q(\bar w)<q(w)<q(\bar w)+\eta\), where \(\psi(s)\coloneqq\varrho s^{1-\theta}\).
			\end{definition}
			
			As will be detailed in \cref{thm:ECS:global}, global convergence is established when the model \(\func{\M}{\R^{nN}\times\R^{nN}}{\Rinf}\) defined as \(\M(\bm w,\bm x)\coloneqq\Phi(\bm w)+\D*(\bm w,\bm x)\) has the KL property.
			The next assumption provides easily verifiable requirements in terms of $f_i$, $h_i$, and \(\varphi\).
			
			\begin{assumption}[global convergence requirements]\label{ass:KL}%
				In problem \eqref{eq:P},
				\bl{\begin{enumeratass}
				\item\label{ass:Kl:C2}%
				{\color{black}%
					for $i\in[N]$, \(f_i,\h_i\in\C^2(\R^n)\) (hence $C=\R^n$) with \(\nabla^2\h_i\succ 0\);%
				}%
				\item\label{ass:Kl:C1}%
				{\color{black}%
					\(\varphi\) has the KL property with exponent \(\theta\in(0,1)\) (\eg when \(f_i\) and \(g\) are semialgebraic) and is level bounded.%
				}%
				\end{enumeratass}}%
			\end{assumption}
			
			\begin{theorem}[global and linear convergence with essentially cyclic rule \eqref{eq:ECS}]\label{thm:ECS:global}%
				Suppose that \cref{ass:basic,ass:KL} are satisfied and that \(g\) is convex.
				Then, the following hold for the iterates generated by \cref{alg:Finito} with an essentially cyclic rule \eqref{eq:ECS}:
				\begin{enumerate}
				\item\label{thm:ECS:global:z}%
					\(\seq{z^k}\) converges to a stationary point \(z^\star\) for \(\varphi\).
				\item
					If \(\theta>\nicefrac12\), then there exists \(c>0\) such that
					\(
						\varphi(z^k)-\varphi(z^\star)
					{}\leq{}
						ck^{-\frac{1}{2\theta-1}}
					\)
					hold for all \(k\in\N\).%
				\item\label{thm:ECS:global:linear}%
					If \(\theta\in(0,\nicefrac12]\), then \(\seq{z^k}\) and \(\seq{\varphi(z^k)}\) converge at \(R\)-linear rate.
				\end{enumerate}
				\begin{proof}
					\bl{%
						Notice that \cref{ass:Kl:C2} along with level boundedness of \(\varphi\) in \cref{ass:Kl:C1} ensures that the requirement in \cref{ass:ECS:subseq:local} is satisfied.
						By the the first claim of \cref{thm:ECS:subseq}, \(\seq{\D**_i(s_i^k,\nabla\h*_i(z^k))}\) converges to zero, and thus we may invoke
					}%
					\cref{thm:subseq} to conclude that the set \(\omega\) of limit points of \(\seq{z^k}\) is nonempty, compact, connected, and made of stationary points for \(\varphi\), with
					\(
						\varphi
					{}\equiv{}
						\varphi_\star
					{}\coloneqq{}
						\lim_{k\to\infty}\FBE(\bm x^k)
					\)
					on \(\omega\).
					If \(\FBE(\bm x^k)=\varphi_\star\) holds for some \(k\in\N\), then it follows from \cref{thm:Igeq} that \(\seq{\bm x^k}\) is asymptotically constant, and the assertion holds trivially.
					In what follows we thus assume that \(\FBE(\bm x^k)>\varphi_\star\) holds for all \(k\).
					The assumptions together with \cref{thm:equivP:KL} ensure that \(\Phi\) enjoys the KL property with exponent \(\theta\).
					Since \(\h*\) is locally strongly convex, we may invoke \cite[Lem. 5.1]{yu2019deducing} to infer that the function \(\func{\M}{\R^{nN}\times\R^{nN}}{\Rinf}\) defined as \(\M(\bm w,\bm x)=\Phi(\bm w)+\D*(\bm w,\bm x)\) has the KL property with exponent \(\vartheta\coloneqq \max\set{\theta,\nicefrac12}\) at every point of the compact set \(\bm\Omega\coloneqq\set{(\bm z^\star, \bm z^\star)}[\bm z^\star\in\bm\omega]\),
					\bl{%
						where \(\bm \omega \coloneqq \set{\bm z = (z, \ldots, z)}[z\in\omega]\).%
					}%
					\footnote{%
						Consistently with the locality of the KL property and the compactness of \(\bm\omega\), the global strong convexity requirement in \cite[Lem. 5.1]{yu2019deducing} can clearly be replaced by local strong convexity.
						Similarly, if \(\Phi\) is a KL function with exponent \(\theta\), then it is trivially a KL function with exponent \(\vartheta\), thus complying with the requirement in the reference.%
					}
					Notice that
					\(
						\FBE(\bm x^k)
					{}={}
						\M(\bm u^k,\bm x^k)
					\)
					and
					\(
						\M(\bm z^\star,\bm z^\star)
					{}={}
						\FBE(\bm z^\star)
					{}={}
						\varphi_\star
					\)
					hold for every \(k\in\N\) and \(\bm z^\star\in\bm\omega\) (cf. \cref{thm:ECS:subseq}), and that
					\(
						\partial\M(\bm w,\bm x)
					{}={}
						\left(\partial\Phi(\bm w) + \nabla\h*(\bm w)-\nabla\h*(\bm x), \nabla^2\h*(\bm x)(\bm x-\bm w)\right)
					\).
					\bl{%
						By \cref{thm:OC} we have \(\nabla\h*(\bm x^k)-\nabla\h*(\bm u^k)\in\partial\Phi(\bm u^k)\), which in turn implies
					}%
					\begin{equation}\label{eq:ECS:distpartial}
						\dist\bigl(0,\partial\M(\bm u^k,\bm x^k)\bigr)
					{}\leq{}
						\|\nabla^2\h*(\bm x^k)\|
						\|\bm x^k-\bm u^k\|
					{}\leq{}
						C\|\bm x^k-\bm u^k\|,
					\end{equation}
					where \(C=\sup_k\|\nabla^2\h*(\bm x^k)\|\) is finite due to boundedness of \(\seq{\bm x^k}\) and continuity of \(\nabla^2\h*\).
					Let \(\psi(t)\coloneqq \rho t^{1-\vartheta}\) be a desingularizing function for \(\M\) on \(\bm\Omega\) \cite[Lem. 1(ii)]{attouch2009convergence}, namely such that
					\[
						\psi'\bigl(\M(\bm w,\bm x)-\varphi_\star\bigr)
						\dist\bigl(0,\partial\M(\bm w,\bm x)\bigr)
					{}\geq{}
						1
					\]
					holds for some \(\varepsilon>0\) and all \((\bm w,\bm x)\) \(\varepsilon\)-close to \(\bm\Omega\) such that \(0<\M(\bm w,\bm x)-\varphi_\star<\varepsilon\).
					Since \(\M(\bm u^k,\bm x^k)=\FBE(\bm x^k)\searrow\varphi_\star\) (cf. \cref{thm:decrease}) and \(\seq{\bm u^k,\bm x^k}\) is bounded and accumulates on \(\bm\Omega\), by discarding early iterates we may assume that the inequality above holds for \((\bm w,\bm x)=(\bm u^k,\bm x^k)\), \(k\in\N\), which combined with \eqref{eq:ECS:distpartial} results in%
					\begin{equation}\label{eq:KL:model}
						\rho^{-1}(1-\theta)^{-1}\bigl(\FBE(\bm x^k)-\varphi_\star\bigr)^\theta
					{}={}
						\psi'\bigl(\FBE(\bm x^k)-\varphi_\star\bigr)^{-1}
					{}\leq{}
						C\|\bm x^k-\bm u^k\|.
					\end{equation}
					Let
					\(
						\Delta_k
					{}\coloneqq{}
						\psi(\FBE(\bm x^k)-\varphi_\star)
					\),
					so that
					\(
						\FBE(\bm x^k)-\varphi_\star
					{}={}
						(\Delta_k/\rho)^{\nicefrac{1}{1-\theta}}
					\).
					By concavity of $\psi$ we have
					\begin{align*}
						\Delta_{T(\nu+1)}-\Delta_{T\nu}
					{}\leq{} &
						\psi'(\FBE(\bm x^{T\nu})-\varphi_\star)
						\bigl(\FBE(\bm x^{T(\nu+1)})-\FBE(\bm x^{T\nu})\bigr)
					{}\overrel*[\leq]{\eqref{eq:KL:model}}{}
						~
						\frac{\FBE(\bm x^{T(\nu+1)})-\FBE(\bm x^{T\nu})}{C\|\bm u^{T\nu}-\bm x^{T\nu}\|}
					\\
					\numberthis\label{eq:ECS:KL:descent}
						\dueto{\eqref{eq:ECS:SD}}
					{}\leq{} &
						-c\|\bm u^{T\nu}-\bm x^{T\nu}\|,
					\end{align*}
					for some constant \(c>0\).
					\bl{%
						As argued in the proof of \cref{thm:ECS:subseq}, by suitably shifting, we conclude that for all \(k\in\N\)
					}%
					\begin{equation}\label{eq:ECS:Delta:shifted}
						\Delta_{k+T}-\Delta_k
					{}\leq{}
						-c\|\bm u^k-\bm x^k\|.
					\end{equation}
					\bl{%
						The rest of the proof is standard (see \cite[Thm. 2]{attouch2009convergence}), and is provided for completeness.
					}%
					Since \(\Delta_k\geq0\), by telescoping we conclude that $\seq{\|\bm u^k-\bm x^k\|}$ has finite sum, and since $\|\bm x^{k+1}-\bm x^k\|\leq \|\bm u^k-\bm x^k\|$ for all $k\in\N$, $\seq{\bm x^k}$ has finite length.
					Therefore, $\seq{\bm x^k}$ and $\seq{\bm u^k}$ converge to the same stationary point, be it \(\bm z^\star\), owing to \cref{thm:ECS:subseq}.
			
					We now show the convergence rates.
					It follows from \eqref{eq:KL:model} that
					\[
						C\rho^{\nicefrac{1}{1-\vartheta}}(1-\vartheta)
						\|\bm x^k-\bm u^k\|
					{}\geq{}
						\rho^{\nicefrac{\vartheta}{1-\vartheta}}
						\big(\FBE(\bm x^{k})-\varphi_\star\big)^\vartheta
					{}={}
						\Delta_k^{\nicefrac{\vartheta}{1-\vartheta}}.
					\]
					Combined with \eqref{eq:ECS:Delta:shifted}, it results in
					\(
						c_1\Delta_k^{\nicefrac{\vartheta}{1-\vartheta}}
					{}\leq{}
						\Delta_k-\Delta_{k+T}
					\)
					for some \(c_1>0\).
					We may now invoke \cref{thm:seqsublinear} to infer that, for every \(t\in[T]\), \(\seq{\Delta_{t+\nu T}}[\nu\in\N]\) converges \(Q\)-linearly (to 0) if \(\theta\leq\nicefrac12\) (which corresponds to \(\vartheta=\nicefrac12\)), and \(\Delta_{t+\nu T}\leq c_3\nu^{-\frac{1-\theta}{2\theta-1}}\) for some \(c_3>0\) otherwise (that is, if \(\vartheta=\theta>\nicefrac12\)).
					Note that the former case implies that \(\seq{\Delta_k}\) converges \(R\)-linearly, whereas the latter implies that \(\Delta_k\leq c_4k^{-\frac{1-\theta}{2\theta-1}}\) holds for every \(k\) and some \(c_4\geq c_3\).
					Recalling that
					\(
						\FBE(\bm x^k)-\varphi_\star
					{}={}
						(\Delta_k/\rho)^{\nicefrac{1}{1-\theta}}
					\),
					the claimed rates of for the cost function follow from \cref{thm:leq}.
					Similarly, when \(\seq{\Delta_k}\) converges \(Q\)-linearly then so does \(\seq{\|\bm x^k-\bm u^k\|}\), as it follows from \eqref{eq:ECS:Delta:shifted}, and in turn so does \(\seq{\|\bm x^k-\bm x^{k+1}\|}\) owing to the inequality \(\|\bm x^k-\bm x^{k+1}\|\leq\|\bm x^k-\bm u^k\|\).
					These two facts imply that \(\seq{\bm x^k}\) and \(\seq{\bm u^k=(z^k,\dots,z^k)}\) are \(R\)-linearly convergent.
				\end{proof}
			\end{theorem}

		\subsection{Low-memory variant}\label{sec:LBfinito}%

			We now analyze \cref{alg:LM}, which, as shown next, is simply a particular implementation of \cref{alg:Finito}. 
			
			\begin{lemma}[\cref{alg:LM} as an instance of \cref{alg:Finito}]\label{thm:equivLM}%
				As long as the same parameters are chosen in \cref{alg:LM,alg:Finito}, to any sequence \(\seq{\s_\LM^k,z_\LM^k,\mathcal I_\LM^{k+1}}\) generated by \cref{alg:LM} there corresponds an identical sequence \(\seq{\s^k,z^k,\mathcal I^{k+1}}\) generated by \cref{alg:Finito}.
				Moreover, the indices \(\seq{\mathcal I_\LM^{k+1}}\) comply with the essentially cyclic rule \eqref{eq:ECS} with \(T=N\).%
				\begin{proof}
					\bl{See \cref{proof:thm:equivLM}.}
				\end{proof}
			\end{lemma}
			
			As a consequence of \cref{thm:equivLM}, \cref{alg:LM} inherits all the convergence results of \cref{Sec:ECS}.
			In addition, here, the convexity requirement of $g$ in \cref{thm:ECS:subseq,thm:ECS:global} can be lifted thanks to the periodic full sampling of the indices.
			
			\begin{theorem}[subsequential convergence of \cref{alg:LM}]\label{thm:LM:subseq}%
				Suppose that \cref{ass:basic} holds and let $\omega$ be the \bl{limit set} of the sequence $\seq{\tilde z^k}$ generated by \cref{alg:LM}.
				Then,%
				\begin{enumerate}
				\item\label{thm:LM:cost}%
					the sequence \(\seq{\varphi(\tilde z^k)}\) converges to the finite value \(\varphi_\star\leq\varphi(x^{\rm init})\) as in \cref{thm:decrease};
				\item\label{thm:LM:const}%
					if \cref{ass:d} is satisfied, then \(\varphi\equiv\varphi_\star\) on \(\omega\);
				\item\label{thm:LM:stationary}%
					if \(C=\R^n\), then \(0\in\hat\partial\varphi(z^\star)\) for every \(z^\star\in\omega\).
				\end{enumerate}
				\begin{proof}
					As shown in \cref{thm:equivLM}, \cref{alg:LM} coincides with \cref{alg:Finito} with an essentially cyclic sampling rule \eqref{eq:ECS}, and there exists an indexing subsequence \(\seq{k_r}[r\in\N]\) with \(0<k_{r+1}-k_r\leq N+1\) such that $\mathcal K^{k_r}=\emptyset$.
					Then, the \(\tilde z\)-update rule (cf. \cref{state:LM:tz_full,state:LM:tz}) yields
					\begin{equation}\label{eq:LM:tz}
						z^{k_r}
					{}={}
						\tilde z^{k_r}
					{}={}
						\tilde z^{k_r+1}
					{}=\dots={}
						\tilde z^{k_{r+1}-1}
					\quad
						\forall r\in\N.
					\end{equation}
					We have
					\begin{align*}
						\FBE(\bm x^{k_{r}+1})
					{}\leq{} &
						\FBE(\bm x^{k_{r}})
						{}-{}
						\sum*_{i\in\mathcal I^{k_{r}+1}}{
							\D*_i(u^{k_{r}},x_i^{k_{r}})
						}
					&
						\dueto{\cref{thm:Igeq}}
					\\
					\numberthis\label{eq:LM:SD}
					{}={} &
						\FBE(\bm x^{k_{r}})
						{}-{}
						\D*(\bm u^{k_{r}},\bm x^{k_{r}})
					&
						\dueto{\(\mathcal I^{k_{r}+1}=[N]\)}
					\\
					{}\leq{} &
						\FBE(\bm x^{k_{r-1}+1})
						{}-{}
						\D*(\bm u^{k_{r}},\bm x^{k_{r}}),
					&
						\dueto{\cref{thm:Igeq}, \(k_{r}\geq k_{r-1}+1\)}
					\end{align*}
					holding for every \(r\in\N\).
					By telescoping and by using the fact that \(\FBE\geq\min\varphi>-\infty\), it follows that \(\seq{\D*(\bm u^{k_{r}},\bm x^{k_{r}})}[r\in\N]\) has finite sum and in particular vanishes.
					Since \(z^{k_r}=\tilde z^{k_r}\),
					\[
						\varphi_\star
					{}\xleftarrow[r\to\infty]{\text{\ref{thm:decrease}}}{}
						\FBE(\bm x^{k_r})
					{}={}
						\varphi(z^{k_r})
						{}+{}
						\D*(\bm u^{k_r},\bm x^{k_r})
					{}={}
						\varphi(\tilde z^{k_r})
						{}+{}
						\D*(\bm u^{k_r},\bm x^{k_r}),
					\]
					whence assertion \ref{thm:LM:cost} follows.
					Assertions \ref{thm:LM:const} and \ref{thm:LM:stationary} follow by patterning the arguments of \cref{thm:constant,thm:stationary}.
				\end{proof}
			\end{theorem}
			
			In the next theorem global convergence results are provided under \cref{ass:KL} in a fully nonconvex setting.
			Moreover, in the spirit of \cref{thm:ECS:global}, linear and sublinear convergence rates are obtained according to the KL exponent.
			
			\begin{theorem}[global and linear convergence of \cref{alg:LM}]\label{thm:LM:global}%
				Suppose that \cref{ass:basic,ass:KL} \bl{hold}.
				Then, the following hold for the iterates generated by \cref{alg:LM}:%
				\begin{enumerate}
				\item\label{thm:LM:global:z}%
					\(\seq{\tilde z^k}\) converges to a stationary point \(x^\star\) for \(\varphi\).
				\item
					If \(\theta>\nicefrac12\), then there exists \(c>0\) such that
					\(
						\varphi(\tilde z^k)-\varphi(x^\star)
					{}\leq{}
						ck^{-\frac{1}{2\theta-1}}
					\)
					hold for all \(r\in\N\).
				\item\label{thm:LM:global:linear}%
					If \(\theta\in(0,\nicefrac12]\), then \(\seq{\tilde z^k}\) and \(\seq{\varphi(\tilde z^k)}\) converge at \(R\)-linear rate.
				\end{enumerate}
				\begin{proof}
					Notice that \cref{ass:KL} entails local strong convexity and Lipschitz differentiability of each \(\h_i\).
					Thus, as discussed in the proof of \cref{thm:ECS:subseq}, $\h*$ is $\mu_{\h*,\bm\U}$-strongly convex on a convex compact set $\mathcal{\bm \U}$ containing the iterates.
					Let the indexing subsequence \(\seq{k_r}[r\in\N]\) be as in the proof of \cref{thm:LM:subseq}, and observe that \(k_r\leq (N+1)r\) holds for every \(r\in\N\).
					The assertions are established \bl{with} the same arguments as in \cref{thm:ECS:global} with the difference that here we examine the generated iterates at subindices $k_r$.
					That is, \eqref{eq:ECS:KL:descent} is replaced by
					\begin{equation}\label{eq:LM:KL:descent}
						\Delta_{k_{r+1}}-\Delta_{k_r}
					{}\leq{}
						\frac{\FBE(\bm x^{k_{r+1}})-\FBE(\bm x^{k_r})}{C\|\bm u^{k_{r}}-\bm x^{k_{r}}\|} 
					{}\underrel*[\leq]{\ref{thm:Igeq}}{}
						\frac{\FBE(\bm x^{k_{r}+1})-\FBE(\bm x^{k_r})}{C\|\bm u^{k_{r}}-\bm x^{k_{r}}\|}
					{}\leq{}
						-C\mu_{\h*,\bm\U}
						\|\bm u^{k_r}-\bm x^{k_r}\|,
					\end{equation}
					where the  last inequality follows from \eqref{eq:LM:SD} and \cref{thm:hstrconvex}.
					Subsequently,
					by patterning the arguments of \cref{thm:ECS:global}, we obtain that $\seq{\Delta_{k_r}}[r\in\N]$ converges $Q$-linearly if $\theta \leq \nicefrac12$, and  $\Delta_{k_r}\leq c r^{-\frac{1-\theta}{2\theta-1}}$ for some $c>0$, if $\theta>\nicefrac12$.
					The claims follow by noting that $\tilde z^k= z^{k_r} = u^{k_r}$ for all $k$ satisfying $k_r \leq k <k_{r+1}$, cf. \eqref{eq:LM:tz}, and arguing as in the last part of \cref{thm:ECS:global}.
				\end{proof}
			\end{theorem}

	\section{Application to phase retrieval and numerical simulations}\label{sec:simulations}%
		In this section we study two examples related to the phase retrieval problem, which consists of recovering a signal based on \emph{intensity measurements}, and arises in many important applications including X-ray crystallography, speech processing, electron microscopy, astronomy, and optical imaging; see, e.g., \cite{candes2015phase,luke2002optical,shechtman2015phase,sun2018geometric}.
		Here, we consider phase retrieval problems with real-valued data, that is, given $a_i\in\R^n\setminus\set0$ and scalars $b_i\in\R_+$, $i\in[N]$, the goal is to find $x\in\R^n$ such that
		\begin{equation}\label{eq:phaseR}
			b_i
		{}\approx{}
			\innprod{a_i}{x}^2, \quad i\in[N],
		\end{equation}
		accounting for the fact that in real-world applications the recorded intensities are likely corrupted by noise, and may involve outliers due to measurement errors.
		To tackle such problems, we consider the following \emph{sparse phase retrieval} formulation:%
		\begin{equation}\label{eq:regPR}
		\textstyle
			\minimize_{x\in\R^n}{
				\tfrac1N
				\sum_{i=1}^N\loss(b_i,\innprod{a_i}{x}^2)+g(x)
			},
		\end{equation}
		where $\loss$ is a loss function, and $g$ is \bl{a} sparsity inducing function (e.g., $l_1$- or $l_0$-norm).
		In particular, we study the case of squared loss $\loss(y,z)=\tfrac{1}{4}(y-z)^2$ \cite{candes2015phase,sun2018geometric}, and Poisson loss $\loss(y,z)=z - y\log(z)$ \cite{chen2015solving,zhang2016provable}, suitable when measurements follow the Poisson model ($b_i\approx\poisson(\innprod{a_i}{x}^2)$).
		Other formulations with $l_1$-loss have been studied in the literature \cite{duchi2019solving,davis2020nonsmooth}.

		\subsection{Sparse phase retrieval with squared loss}\label{sec:Quartic}%
			Consider the nonconvex minimization \eqref{eq:regPR} with squared loss $\loss(y,z)=\tfrac{1}{4}(y-z)^2$, and either $g=\lambda\|{}\cdot{}\|_1$, $\lambda \geq 0$, or $g=\indicator_{\mathbb B_\kappa}$, where $\mathbb{B}_\kappa$ is the $l_0$-norm ball of radius $\kappa$.
			This problem is written in the form of \eqref{eq:P} with
			\begin{equation}\label{eq:fiSQIP}
				f_i(x)
			{}={}
				\tfrac14
				(\innprod*{a_i}{x}^2-b_i)^2
			\quad\text{and}\quad
				\h_i(x)
			{}={}
				\tfrac14\|x\|^4
				{}+{}
				\tfrac12\|x\|^2.
			\end{equation}
			The next lemma is a simple adaptation from \cite{bolte2018first} for finding the smoothness moduli of $f_i$ relative to the Legendre kernel $\h_i$, and for computing the solutions to the subproblem \eqref{eq:innerprox}.
			For $l_1$-regularization, the inner subproblems amount to computations involving the soft-thresholding operator, whereas in the case of $l_0$-norm ball they amount to computing projections onto $\mathbb B_\kappa$, that is, setting to zero $n-\kappa$ \bl{elements among the} smallest in absolute value.
			
			\begin{lemma}[{{\cite[Lem. 5.1, Props. 5.1 and 5.2]{bolte2018first}}}]\label{lem:plip}
				Let $f_i$ and $\h_i$ be as in \eqref{eq:fiSQIP}.
				Then, $f_i$ is $L_{f_i}$-smooth relative to $\h_i$ with $L_{f_i}=(3\|a_i\|^4+\|a_i\|^2 |b_i|)$.
				Moreover, denoting
				\(
					\bar\gamma
				{}={}
					(\sum_{i=1}^N\nicefrac{1}{\gamma_i})^{-1}
				\),
				for any
				\(
					y(s)
				{}\in{}
					\prox_{\bar\gamma g}(\bar\gamma s)
				\)
				the operator \(T\) as defined in \eqref{eq:innerprox} may be computed as follows:
				\begin{enumerate}
				\item
					If $g=\lambda\|{}\cdot{}\|_1$, $\lambda\geq0$, then \(T(s)=t^\star y(s)\), where \(t^\star\) is the real positive root of the equation
					\(
						\|y(s)\|^2t^3+t-1=0
					\).\footnote{\label{fn1}%
						Nonnegative real roots of the cubic equation $t^3+pt+q=0$ for some $p>0$ and $q\leq 0$ are given by Cardano's formula
						\(
							t^\star
						{}={}
							(c-\nicefrac q2)^{\nicefrac13}
							{}-{}
							\big(c+\nicefrac q2)^{\nicefrac13}
						\),
						where
						\(
							c
						{}={}
							(\nicefrac{q^2}{4}+\nicefrac{p^3}{27})^{\nicefrac12}
						\);
						see, e.g., \cite{spiegel1999mathematical}.%
					}
				\item
					If \(g=\indicator_{\mathbb{B}_\kappa}(x)\), then \(T(s)\ni-t^\star\|y(s)\|^{-1}y(s)\), where \(t^\star\) is the real nonnegative root of \(t^3+t- \|y(s)\|=0\) (see \cref{fn1}).
				\end{enumerate}
			\end{lemma}

		\subsection{Sparse phase retrieval with Poisson loss}\label{sec:Poission}%
			We now assume that the recorded intensities follow the Poisson model ($b_i\sim\poisson(\innprod{a_i}{x}^2)$).
			In this setting we adapt the Poisson loss $\loss(y,z)=z - y\log(z)$ and consider the  regularized problem \eqref{eq:regPR} with $g=\lambda\|{}\cdot{}\|_1$, $\lambda\geq 0$.
			This problem may be written in the form of \eqref{eq:P} by setting
			\begin{equation}\label{eq:f:PhR:Poisson}
				f_i(x)
			{}={}
				-b_i\log(\innprod{a_i}{x}^2)
				{}+{}
				\innprod{a_i}{x}^2
			\quad\text{and}\quad
				\h_i(x)
			{}={}
				\|a_i\|^2\|x\|^2
				{}-{}
				\smash{2b_i\sum_{j=1}^n\log(x_j).}
			\end{equation}
			As shown next, the nonconvex function $f_i$ is smooth relative to $\h_i$, and the operator \(T\) as in \eqref{eq:innerprox} is easily computable.

			\begin{lemma}
				Let \(f_i\) and \(\h_i\) be as in \eqref{eq:f:PhR:Poisson}, with $a_i\in\R^n_{+}\setminus\set0$ and $b=(b_1,\ldots,b_N)\in\R^n_{+}\setminus\set0$. 
				Then, \(\h_i\) is a \(2\|a_i\|^2\)-strongly convex Legendre kernel (with \(\dom\h_i=\R_{++}^n\)), and \(f_i\) is \(L_{f_i}\)-smooth relative to \(\h_i\) with \(L_{f_i}=1\).
				Moreover, denoting
					\(
					c_a
				{}={}
					\sum_{i=1}^N\tfrac{4}{\gamma_i}\|a_i\|^2
				\)
				and
				\(
					c_b
				{}={}
					\sum_{i=1}^N\tfrac{4}{\gamma_i}b_i
				\),
				the operator \(T\) as defined in \eqref{eq:innerprox} with \(g=\lambda\|{}\cdot{}\|_1\), $\lambda\geq 0$, is given by
				\begin{equation}\label{eq:prox:PhR:Poisson:ell1}
					T(s)
				{}={}
					(w_1,\ldots,w_M),
				\quad\text{with}\quad
					w_j
				{}={}
					\tfrac{1}{c_a}\Bigl(%
						s_{j}
						{}-{}
						\lambda
						{}+{}
						\bigl(
						(s_{j}-\lambda)^{2}
						{}+{}
						c_a c_b\bigr)^{\nicefrac{1}{2}}
					\Bigr).
				\end{equation}
				\begin{proof}
					To avoid clutter, we drop the \bl{subscripts} $i$.
					The assertion on \(h\) is of immediate verification.
					Since both \(f\) and \(h\) are \(\C^2\) on \(\interior\dom h\), once we show that \(\nabla^2h-\nabla^2f\succeq0\) on \(\interior\dom h\) the claim will follow. 
					From direct computations,
					\(
						\nabla^2h(x)
					{}={}
					   2\|a\|^2
						{}+{}
						2b\diag(x_1^{-2},\dots,x_n^{-2})
						\)
					and
					\(
						\nabla^2f(x)
					{}={}
						2
						\bigl(
							1+\tfrac{b}{\innprod ax^2}
						\bigr)
						a\trans a
					\).
					In particular,
					\(
						M(x)
					{}\coloneqq{}
						\nabla^2h(x)-\nabla^2f(x)
					{}={}
						2b\bigl(
							\diag(x_1^{-2},\dots,x_n^{-2})
							{}-{}
							\tfrac{1}{\innprod ax^2}
							a\trans a
						\bigr)
						{}+{}
						2\|a\|^2
						{}-{}
						2a\trans a
					\).
					For every \(y\in\R^n\) it holds that (here $a_k$ is the $k$-th coordinate of $a$)
					\[
						\innprod*{y}{M(x)y}
					{}\geq{}
						2b\sum_{j=1}^n{
							\frac{y_j^2}{x_j^2}
						}
						{}-{}
						2b\frac{\innprod ay^{2}}{\innprod ax^{2}}
					{}\geq{}
						2b\sum_{j=1}^{n}\frac{y_{j}^{2}}{x_{j}^{2}}
						{}-{}
						2b\sum_{j=1}^{n}\frac{a_{j}x_{j}}{\innprod ax}\frac{y_{j}^{2}}{x_{j}^{2}}
					{}\geq{}
						0,
					\]
					where the first inequality follows by the Cauchy-Schwarz inequality,
					\bl{%
						the second one from Jensen's inequality
						\(
							\langle a, y \rangle^2 = \sum_{j=1}^n \big(a_j x_j (\tfrac{y_j}{x_j})\big)^2 \leq \sum_{i=1}^n a_i x_i \sum_{j=1}^n a_j x_j (\tfrac{y_j}{x_j})^2
						\),
					}%
					and \bl{the} third one from the fact that \(\sum_i\alpha_i\sum_j\beta_j\geq\sum_i\alpha_i\beta_i\) for every \(\alpha_i,\beta_i\geq0\), \(i\in[n]\).
					The closed-form solution for the proximal mapping $T(s)$ follows directly from its first\bl-order optimality conditions.
				\end{proof}
			\end{lemma}

		\subsection{Experimental setup}\label{sec:PhR}%
			We test \cref{alg:Finito} with sampling rules \eqref{eq:RS} (using single-index selection with uniform sampling), \eqref{eq:shuffled}, \eqref{eq:cyclic}, and the low-memory \cref{alg:LM} with a cyclic inner loop (corresponding to \(\mathcal I^{k+1}=[N]\) if \(\mod(k, N+1)=0\), and \(\mathcal I^{k+1}=\set{\mod(k,N+1)}\) otherwise).
			\bl{%
				We also consider the full mirror descent algorithm (MD) for nonconvex problems under the relative smoothness assumption \cite{bolte2018first}, and SMD, its stochastic extension \cite{davis2018stochastic,hanzely2018fastest}.
			}%
			The incremental method PLIAG \cite{zhang2021proximal} was not tested, as the problem setting does not comply with the requirements therein; cf. \cite[\bl{Assump. 8}]{zhang2021proximal}.
			\bl{%
				For the problem of \cref{sec:Quartic} with \(g=\indicator_{\mathbb B_\kappa}\), the (shuffled) cyclic rules in \cref{alg:Finito} do not comply with \cref{thm:ECS:subseq} (since \(g\) is nonconvex), but are however provided as empirical evidence.
			}%
			
			\paragraph{Parameters selection}
				For \cref{alg:Finito,alg:LM}, we always use $\gamma_i = \nicefrac{0.99N}{L_{f_i}}$.
				For SMD we used the popular square-summable stepsize $\gamma_k=\nicefrac{\alpha}{(L_f k)}$, where $k$ is the iteration counter, $L_f$ is the smoothness modulus of $f=\nicefrac1N\sum_{i=1}^Nf_i$ relative to a suitable Bregman kernel $h$, and $\alpha>0$ is tuned for performance.
				In particular, for SMD in the problems described above, $L_f = \sum_{i=1}^N \frac1N{L_{f_i}}$ and $h(x)=\frac14\|x\|^4 + \frac12\|x\|^2$ for simulations related to \cref{sec:Quartic}, and $h(x)= \frac1N\sum_{i=1}^N\|a_i\|^2\|x\|^2 - \frac2N\sum_{i=1}^Nb_i\sum_{j=1}^n\log(x_j)$ for those related to \cref{sec:Poission} are used.
			
			\paragraph{Optimality criteria}
				As a measure of suboptimality, we consider
				\begin{equation}\label{eq:residual_z}
				\textstyle
					\Op(z^k)
				{}\coloneqq{}
					\|z^k{}-{}v^k\|
				\quad\text{for some}\quad
					v^k\in T(\textstyle\sum_{i=1}^N\nabla\h*_i(z^k)),
				\end{equation}
				since it satisfies
				\begin{align}
				\label{eq:distvarphi}
					\tfrac1N\dist(0,\hat\partial\varphi(v^k))
				{}\overrel*{\ref{thm:equivP:partial}}{} &
				\textstyle
					\inf_{\bm w\in  \hat \partial\Phi(\bm v^k)}\tfrac1N\|\sum_{i=1}^N w_i\|
				{}\overrel*[\leq]{\ref{thm:OC}}{}
					\tfrac1N\bigl\|
						\sum_{i=1}^N\bigl(
							\nabla\h*_i(z^k)
							{}-{}
							\nabla\h*_i(v^k)
						\bigr)
					\bigr\|
				\\
				\nonumber
				{}\leq{} &
				\textstyle
					\tfrac1N\sum_{i=1}^N\|\nabla\h*_i(z^k)
					{}-{}
					\nabla\h*_i(v^k)\|
				{}\leq{}
					\eta
					 \|z^k-v^k\|
				{}={}
					\eta
					\Op(z^k),
				\end{align}
				where $\bm v^k=(v_1,\ldots, v_N)$, and $\eta$ is some positive constant that exists by virtue of local Lipschitz continuity of $\h_i$, \cref{thm:Legendre:smooth}, and boundedness of the sequences $\seq{z^k}$ and $\seq{v^k}$.
				Although, in a similar fashion to \eqref{eq:distvarphi}, \(\dist(0,\hat\partial\varphi(z^k))\) may be upper bounded by \(\|\tilde s^k - \sum_i\nabla \h*_i(z^k)\|\) which would  be an equally good estimate of \(\dist(0,\hat\partial\varphi(z^k))\), this quantity is not readily available in other methods such as SMD.
				The introduction of \(v^k\), instead, offers a viable algorithm-independent alternative that only requires access to output variables.\footnote{%
				\bl{%
					Using \emph{almost sure nondegeneracy} of the fixed points of \(\prox^{\h*}_\Phi\) (cf. \cite[Def. 3.5 and Lem. 3.6]{ahookhosh2021bregman} and discussion therein), hence that of \(T\circ\sum_i\nabla\h_i\) by \cref{thm:Legendre:rangeT}, it can be deduced that \(\Op(z^k)\) vanishing is necessary for optimality of the limit point(s) of \(z^k\).
					We however omit the technical details, as this behavior is anyway confirmed in the plots.%
				}\label{footnote:Op}}

			\subsubsection*{Simulations}%
				In the first set of simulations we consider $16 \times 16$ gray-scale images from a digits dataset \cite{friedman2001elements}\footnote{%
					\url{https://web.stanford.edu/~hastie/ElemStatLearn/data.html}.%
				}
				and a QR code dataset \cite{metzler2017coherent}.
				The images are vectorized resulting in the signal $x\in\R^n$ with $n=256$.
				The data matrix $A\in\R^{N\times n}$ ($a_i$ being the $i$-th row) with $N=nd$, $d=5$ is generated following the procedure described in \cite[\S 6.3]{duchi2019solving}.
				Let $M\in\R^{n\times n}$ be a normalized Hadamard matrix.
				We generate $d$ many i.i.d. diagonal sign matrices $S_i$ with diagonal elements in $\set{-1,1}$ selected uniformly at random, and set
				\(
					A
				{}={}
					{[MS_1,\dots,MS_d]}
				\).
				Typically $d\geq 3$ is sufficient for near complete recovery on noiseless data.
				In our simulations, we corrupted a fraction of the measurements $b_i=\innprod*{a_i}{x}^2$ independently by setting $b_i=0$ with probability $p_{\rm c}=\nicefrac{1}{50}$.
				All of the plotted algorithms are initialized using the initialization scheme described in \cite[\S3]{duchi2019solving}.
				
				\begin{figure}[bt]
					\centering
					\includetikz[width=0.75\linewidth]{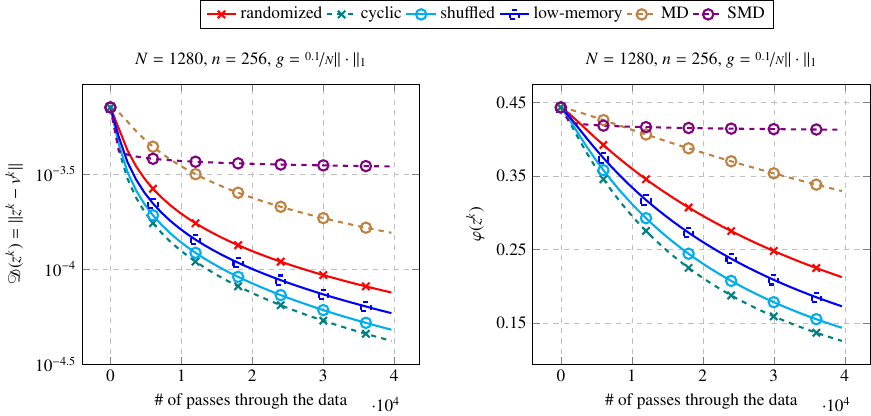}
					\includetikz[width=0.75\linewidth]{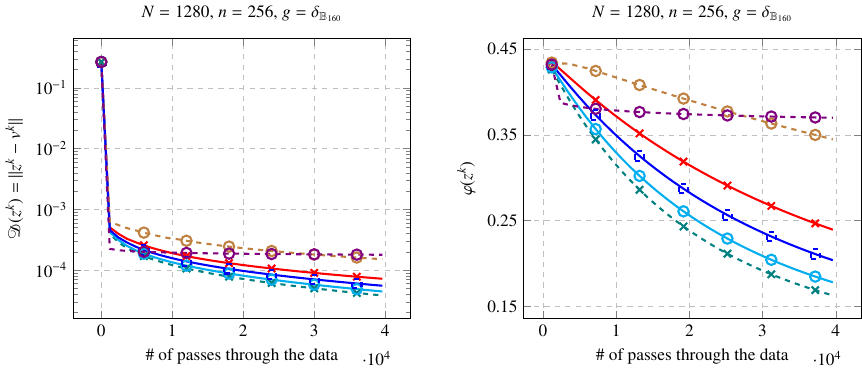}%
					\caption[Representative convergence plots for problem \eqref{eq:regPR} with squared loss on a digits image]{%
						Representative convergence plots for problem \protect\eqref{eq:regPR} with squared loss on a digits image: (first row) $\ell_1$-regularization, (second row) $\ell_0$-norm ball constraint.
						The related plots for the QR code images follow a very similar trend and are therefore omitted.%
					}%
					\label{fig3c}%
				\end{figure}
				
				\begin{figure}
					\centering
					\begin{subfigure}{.2\textwidth}
						\centering
						\includegraphics[scale = 2.1]{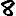}
						
						\includegraphics[scale = 2.1]{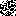}%
						\caption{original image}%
						\label{fig:sub1-QR}
					\end{subfigure}
					\begin{subfigure}{.2\textwidth}
						\centering
						\includegraphics[scale = 2.1]{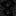}
				
						\includegraphics[scale = 2.1]{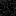}%
						\caption{initialization}%
						\label{fig:sub3-QR}
					\end{subfigure}
					\begin{subfigure}{.2\textwidth}
						\centering
						\includegraphics[scale = 2.1]{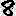}
				
						\includegraphics[scale = 2.1]{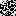}%
						\caption{tolerance $10^{-5}$}%
						\label{fig:sub5-QR}%
					\end{subfigure}
					\begin{subfigure}{.2\textwidth}
						\centering
						\includegraphics[scale = 2.1]{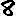}
				
						\includegraphics[scale = 2.1]{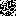}%
						\caption{tolerance $10^{-7}$}%
						\label{fig:sub4-QR}%
					\end{subfigure}
					\caption{%
						Image recovery with corrupted measurements for tolerances $\protect\set*{10^{-5},10^{-7}}$.
						The sparsity parameters $\kappa=160$ and $\kappa=125$ are used for the digit and the QR code, respectively.%
					}%
					\label{fig:Quartic_digits}%
				\end{figure}
				
				For the $l_1$-regularized problem we performed tests with different values of the regularization parameter $\lambda$ and found $\lambda =\nicefrac{0.1}N$ to lead to a visually favorable recovery.
				When $g=\delta_{\mathbb{B}_\kappa}$, we set $\kappa=160$ and $\kappa=125$ for the digit and QR data, respectively.
				The convergence behavior in terms of $\Op(z^k)$ (see \eqref{eq:residual_z}) is plotted in \cref{fig3c} for a representative digit $8$ image.
				With the above described initialization, the algorithms converge to the same cost.
				In our simulations SMD had the slowest performance, and the \bl{cyclic rule} \eqref{eq:cyclic} in \cref{alg:Finito} was observed to consistently outperform all others.
				\bl{%
				The low-memory \cref{alg:LM} has a comparable performance, almost always superior to the randomized variant. 
				}%
				As evident from \cref{fig:Quartic_digits}, despite corrupted measurements a reasonably good recovery is achieved with the $l_0$-norm ball.
				A similar recovery is observed with $l_1$ regularization.

			\begin{figure}[bt]
				\center
				\includetikz[width=0.75\linewidth]{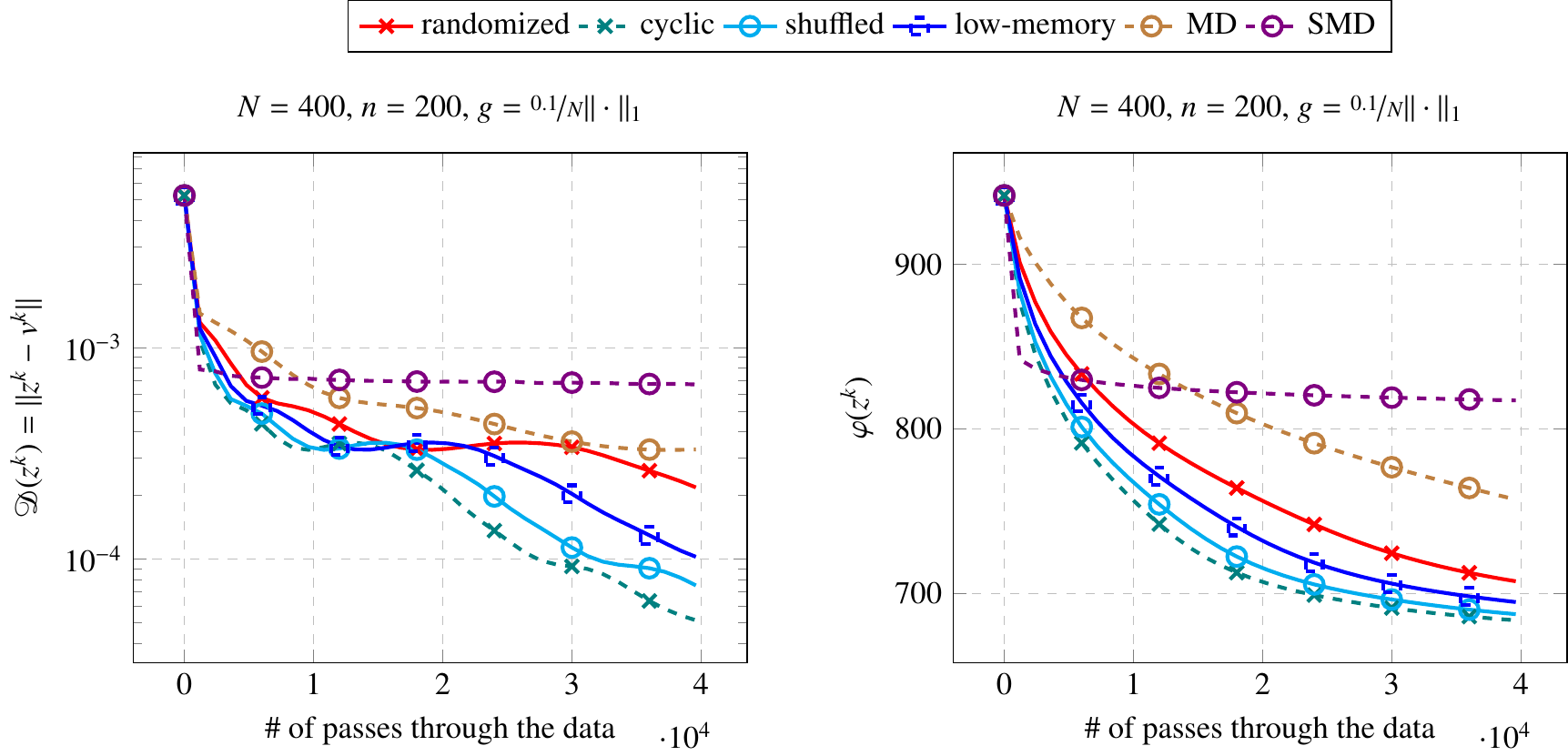}
				\includetikz[width=0.75\linewidth]{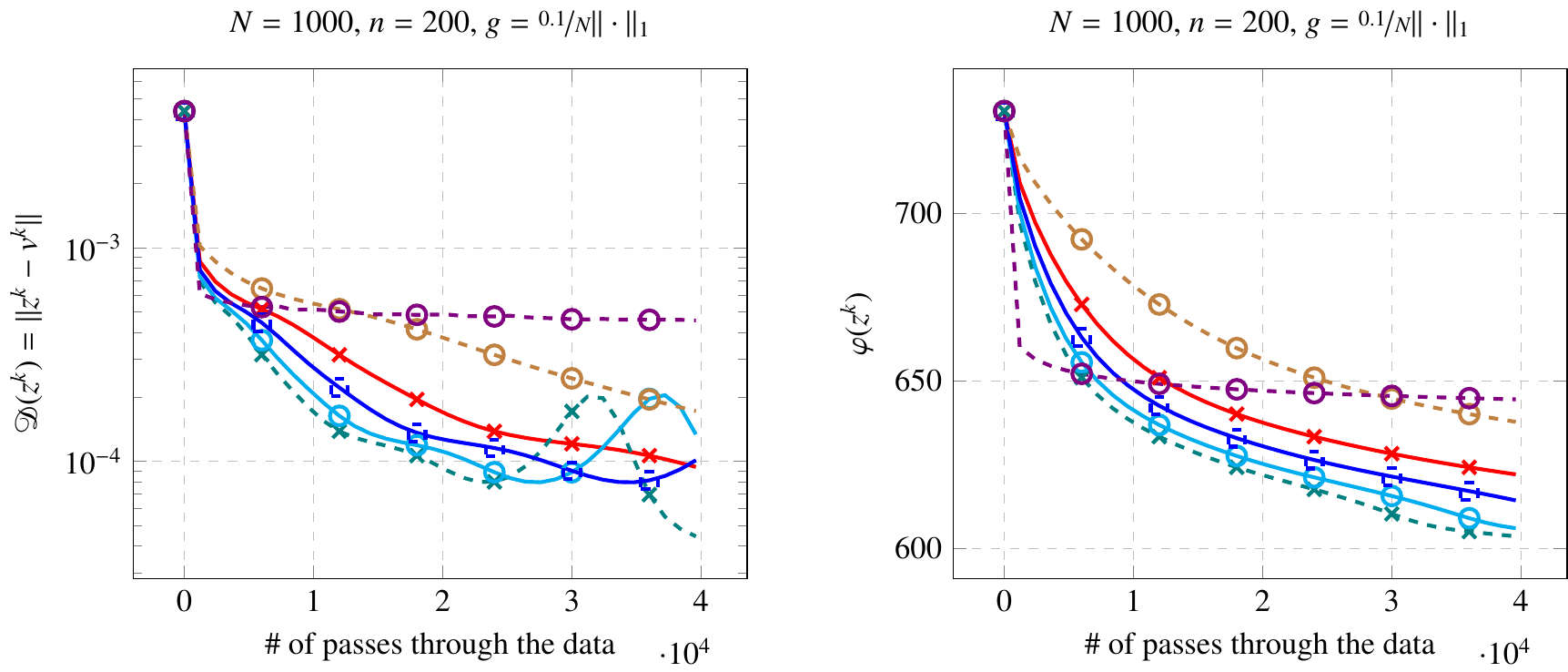}%
				\caption{%
				   Representative convergence plots for the \(l_1\)-regularized problem with Poisson loss.%
				}%
				\label{fig:PRLPoisson}%
			\end{figure}

			In the last set of simulations we consider synthetic data. We generate a standard random Gaussian matrix $A\in\R^{N\times n}$ with  $n=200$, $N\in\set{400,1000}$.
			The data vector $a_i$, $i\in[N]$, is set equal to the absolute value of the $i$-th row of $A$.
			We also drew  a random vector from $\mathcal{N}(0,I_n)$ and set the signal $x$ equal to its absolute value.
			We generated the measurements according to the Poisson model $b_i\sim\poisson(\innprod*{a_i}{x}^2)$, $i\in[N]$, and further corrupted the measurements $b_i$ by setting them equal to the nearest integer to the absolute value of $\|x\|^2\mathcal N(0,1)$ with probability $p_{\rm c}=\nicefrac{1}{10}$.
			All methods were initialized at the same random point.
			We ran simulations with regularization parameter $\lambda\in\set{\nicefrac{0.01}{N}, \nicefrac{0.1}{N}, \nicefrac{1}{N}}$.
			We only report the results for $\lambda=\nicefrac{0.1}{N}$ due to space limitations, nevertheless remarking that for other values similar plots were observed.
			The results are illustrated in \cref{fig:PRLPoisson}.
			Similar to the previous experiments, SMD \bl{performed the} worst, while the best results are observed for \cref{alg:Finito} with cyclic \bl{rule} \eqref{eq:cyclic}.
			The low-memory \cref{alg:LM} usually outperforms the randomized variant of \cref{alg:Finito}.

	\section{Conclusions}\label{sec:conclusions}%
		A Bregman incremental aggregated method was developed that extends Finito/MISO \cite{defazio2014finito,mairal2015incremental} to non-Lipschitz and nonconvex settings.
		The basic algorithm was studied under randomized and essentially cyclic sampling strategies.
		Furthermore, a variant with $O(n)$ memory requirements is developed that is novel even in the Euclidean case.
		A sure descent property established on a Bregman Moreau envelope leads to a surprisingly simple convergence analysis.
		As one particularly interesting result, in the randomized setting linear convergence is established under strong convexity of the cost without requiring convexity of the individual functions $f_i$ or $g$.
		Future research directions include extending the analysis to the framework of the Douglas-Rachford splitting, momentum-type schemes, as well as applications to nonconvex distributed asynchronous optimization.

	\begin{appendix}%
		\section{Auxiliary results}\label{sec:appendix:aux}%
			\begin{lemma}[equivalence between \eqref{eq:bigP} and \eqref{eq:P}]\label{thm:equivP}%
			\bl{%
				Let \(\Phi\) and \(\Delta\) be as in \eqref{eq:bigP}.
				Then:
			}%
				\begin{enumerate}
				\item{\sc cost:}\label{thm:equivP:phi}~
					\(\Phi(\bm x)=\varphi(x)\) if \(\bm x=(x,\ldots,x)\), and \(\Phi(\bm x)=\infty\) otherwise.
				\item{\sc subdifferential:}\label{thm:equivP:partial}~
					\(
						\hat\partial(\Phi+\indicator_{\overline C\times\cdots\times\overline C})(\bm x)
					{}={}
						\set{\bm v=(v_1,\ldots,v_N)\in\R^{nN}}[
							\sum_{i=1}^nv_i
						{}\in{}
							\hat\partial(\varphi+\indicator_{\overline C})(x)
						]
					\)
					if \(\bm x=(x,\cdots,x)\) for some \(x\in\R^n\), and is empty otherwise;
					the same relation still holds if the regular subdifferential \(\hat\partial\) is replaced by the limiting subdifferential \(\partial\).
				\item{\sc KL property:}\label{thm:equivP:KL}~
					\(\varphi\) has the KL property at \(x\) iff so does \(\Phi\) at \(\bm x=(x,\dots,x)\), in which case the desingularizing functions are the same up to a positive scaling.
				\item{\sc stationary points:}\label{thm:equivP:stationary}~
					a point \(\bm x^\star\) is stationary for problem \eqref{eq:bigP} iff \(\bm x^\star=(x^\star,\ldots,x^\star)\) for some \(x^\star\in\R^n\) which is stationary for problem \eqref{eq:P}.
				\item{\sc minimizers:}\label{thm:equivP:argmin}~
					\(\bm x^\star\) is a (local) minimizer of problem \eqref{eq:bigP} iff \(\bm x^\star=(x^\star,\ldots,x^\star)\) for some \(x^\star\in\R^n\) which is a (local) minimizer for problem \eqref{eq:P};
					in fact, \(\inf_{\overline C\times\cdots\times\overline C}\Phi=\inf_{\overline C}\varphi\).
				\item{\sc level boundedness:}\label{thm:equivP:lb}~
					\(\varphi+\indicator_{\overline C}\) is level bounded iff so is \(\Phi+\indicator_{\overline C\times\cdots\times\overline C}\).
				\item{\sc convexity:}\label{thm:equivP:convex}~
					\(\func{\varphi}{\R^n}{\Rinf}\) is convex iff so is \(\func{\Phi}{\R^{Nn}}{\Rinf}\).
				\end{enumerate}
			\bl{%
				\begin{proof}%
					For notational convenience, up to possibly replacing \(g\) with \(g+\indicator_{\overline C}\) we may assume without loss of generality that \(C=\R^n\).
					\begin{proofitemize}
					\item\ref{thm:equivP:phi}~
						Trivial consequence of the fact that \(\dom\Phi\subseteq\Delta\) (the consensus set, cf. \eqref{eq:bigP}).
					\item\ref{thm:equivP:partial}~
						In light of the previous point, having \(\bm x=(x,\ldots,x)\) for some \(x\in\R^n\) is necessary for the nonemptiness of \(\hat\partial\Phi(\bm x)\).
						Let \(\bm x=(x,\ldots,x)\) and \(\bm v\in\hat\partial\Phi(\bm x)\) be fixed.
						Then,
						\begin{equation}\label{eq:liminf}
							0
						{}\leq{}
							\liminf_{\bm x\neq\bm y\to\bm x}{
								\frac{
									\Phi(\bm y)-\Phi(\bm x)
									{}-{}
									\innprod*{\bm v}{\bm y-\bm x}
								}{
									\|\bm y-\bm x\|
								}
							}
						{}={}
							\liminf_{x\neq y\to x}{
								\frac{
									\varphi(y)-\varphi(x)
									{}-{}
									\innprod*{\sum_iv_i}{y-x}
								}{
									\sqrt N\|y-x\|
								}
							}
						\end{equation}
						where the equality comes from the fact that \(\dom\Phi\subseteq\Delta\) together with assertion \ref{thm:equivP:phi}.
						This shows that \(\sum_iv_i\in\hat\partial\varphi(x)\).
						Conversely, let \(u\in\hat\partial\varphi(x)\) and \(\bm v\in\R^{nN}\) be such that \(\sum_iv_i=u\).
						By reading \eqref{eq:liminf} from right to left we obtain that \(\bm v\in\hat\partial\Phi(\bm x)\).
						Having shown the identity of the regular subdifferential, the same claim with the limiting subdifferential follows by definition.
					\item\ref{thm:equivP:KL}~
						Follows from the bounds
						\[
							\tfrac1N
							\dist(0,\partial\varphi(x))^2
						{}={}
							\!\!\!\!
							\inf_{\bm v\in\partial\Phi(\bm x)}{\textstyle
								\tfrac1N
								\|\sum_{i=1}^Nv_i\|^2
							}
						{}\leq{}
							\overbracket*{
							\!\!\!\!
								\inf_{\bm v\in\partial\Phi(\bm x)}{\textstyle
									\sum_{i=1}^N\|v_i\|^2
								}
							}^{
								\dist(0,\partial \Phi(\bm x))^2
							}
						{}\leq{}
							\!\!\!\!
							\inf_{\bm v\in\partial\Phi(\bm x)}{\textstyle
								\|\sum_{i=1}^Nv_i\|^2
							}
						{}={}
							\dist(0,\partial\varphi(x))^2,
						\]
						where the first and last equalities are due to assertion \ref{thm:equivP:partial}.
					\item\ref{thm:equivP:stationary}--\ref{thm:equivP:convex}~
						Directly follow from assertions \ref{thm:equivP:phi} and \ref{thm:equivP:partial}.
					\qedhere
					\end{proofitemize}
				\end{proof}
			}%
			\end{lemma}

			\begin{lemma}\label{thm:reg:env}
				Let \(\bm\U\coloneqq\U_1\times\cdots\times\U_N\) with \(\U_i\subseteq\interior\dom\h_i\) nonempty and convex, \(i\in[N]\).
				Additionally to \cref{ass:basic}, suppose that \(g\) is convex, and \(\h_i\), $i\in[N]$, is \(\ell_{\h_i,\U_i}\)-Lipschitz-differentiable and \(\mu_{\h_i,\U_i}\)-strongly convex on \(\U_i\).
				Then, the following hold for function \(\h*\) as in \eqref{eq:hatH} with \(\gamma_i\in(0,\nicefrac{N}{L_{f_i}})\), \(i\in[N]\):
				\begin{enumerate}
				\item\label{thm:TLip}%
					\(\prox_\Phi^{\h*}\) is Lipschitz continuous on \(\bm\U\).
				\end{enumerate}
				If in addition \(f_i\) and \(h_i\) are twice continuously differentiable on \(\U_i\), \(i\in[N]\), then
				\begin{enumerate}[resume]
				\item\label{thm:C1+}%
					\(\FBE\) is continuously differentiable on \(\bm\U\) with \(\nabla\FBE=\nabla^2\h*\circ(\id-\prox_\Phi^{\h*})\);
				\item\label{thm:distres}%
					\(
						\dist(0,\partial\FBE(\bm x))
					{}={}
						\|\nabla\FBE(\bm x)\|
					{}\leq{}
						C_{\bm\U}
						\|\bm x-\bm z\|
					\)
					for any \(\bm x\in\bm\U\), where \(\bm z=\prox_\Phi^{\h*}(\bm x)\) and
					\(
						C_{\bm\U}
					{}={}
						\max_i\big\{
						\big(1+\tfrac{\gamma_iL_{f_i}}{N}\big)
						\frac{\ell_{\h_i,\U_i}}{\gamma_i}
						\big\}
					\).
				\end{enumerate}
				\begin{proof}
					It follows from \cref{thm:Legendre:smooth,thm:Legendre:strconvex} that \(\h*_i\) is \(\ell_{\h*_i,\U_i}\)-Lipschitz differentiable and \(\mu_{\h*_i,\U_i}\)-strongly convex on \(\U_i\) with
					\(
						\ell_{\h*_i,\U_i}
					{}={}
						\big(1+\tfrac{\gamma_iL_{f_i}}{N}\big)
						\frac{\ell_{\h_i,\U_i}}{\gamma_i}
					\)
					and
					\(
						\mu_{\h*_i,\U_i}
					{}={}
						\big(1-\tfrac{\gamma_iL_{f_i}}{N}\big)
						\frac{\mu_{\h_i,\U_i}}{\gamma_i}
					\).
					\bl{Then,} \(\h*\) is Lipschitz differentiable and strongly convex on \(\bm\U\) (with respective moduli \(\ell_{\h*,\bm\U}=\max_i\ell_{\h*_i,\U_i}\) and \(\mu_{\h*,\bm\U}=\min_i\mu_{\h*_i,\U_i}\)), and therefore so is its conjugate \(\conj{\h*}\) on \(\h*(\bm\U)\) (with respective moduli
					\(
						\ell_{\conj{\h*},\bm\U}
					{}={}
						\mu_{\h*,\bm\U}^{-1}
					\)
					and
					\(
						\mu_{\conj{\h*},\bm\U}
					{}={}
						\ell_{\h*,\bm\U}^{-1}
					\)).
					Notice that convexity of \(g\), this being equivalent to that of \(G\), implies that \(\h*(\bm x)+\Phi(\bm x)= G(\bm x) + \sum_{i=1}^N\tfrac1{\gamma_i}\h_i(x_i)\) is strongly convex on \(\bm\U\).
					We may thus invoke \cite[Thm. 4.2]{kan2012moreau} and \cref{thm:inv1} to conclude that
					\(
						\prox_\Phi^{\h*}
					{}={}
						\partial\conj{(\h*+\Phi)}\circ\nabla\h*
					{}={}
						\nabla\conj{(\h*+\Phi)}\circ\nabla\h*
					\)
					is the composition of Lipschitz-continuous mappings on \(\bm\U\), which shows assertion \ref{thm:TLip}.
					In turn, assertion \ref{thm:C1+} follows from \cite[Cor. 3.1]{kan2012moreau}.
					Finally, \(\ell_{\h*,\bm\U}\)-Lipschitz continuity of \(\nabla\h*\) on \(\bm\U\) entails the bound \(\|\nabla^2\h*\|\leq\ell_{\h*,\bm\U}\) on \(\bm\U\), leading to assertion \ref{thm:distres}.
				\end{proof}
			\end{lemma}
			
			\begin{lemma}\label{thm:seqsublinear}%
				Let \(\seq{\alpha_k}\subset\R_+\) be a sequence, and suppose that there exist \(c>0\) and \(\delta\in[1,\infty)\) such that
				\(
					\alpha_{k+1}^\delta
				{}\leq{}
					c(\alpha_k-\alpha_{k+1})
				\)
				holds for every \(k\in\N\).
				\begin{enumerate}
				\item
					If \(\delta=1\), then \(\seq{\alpha_k}\) is \(Q\)-linearly convergent (to \(0\)).
				\item
					If \(\delta\in(1,\infty)\), then there exists \(c'>0\) such that
					\(
						\alpha_k
					{}\leq{}
						c'k^{-\frac{1}{\delta-1}}
					\)
					holds for all \(k\in\N\).
				\end{enumerate}
			\end{lemma}
	
	\bl{%
		\proofsection{sec:conv}%
			\begin{appendixproof}{thm:RS:strconvex}[ (linear convergence with randomized rule \eqref{eq:RS})]~\par\noindent
				We will use the equivalent BC-reformulation of \cref{alg:BCPP}, through the identities shown in \cref{thm:equivBC}.
				We start by observing that \(x_i^k\in\set*{x^{\rm init},z^k}[k\in\N]\subseteq\U\) holds for any \(k\in\N\) and \(i\in[N]\), as it follows from the \(\bm x\)-update at \cref{state:BCPP:x+} and the fact that \(u_i^k=z^k\), cf. \cref{thm:equivBC:z}.
				Let \(\bm x^\star=(x^\star,\ldots,x^\star)\) be the unique minimizer of \(\Phi\) (cf. \cref{thm:equivP:argmin}).
				As shown in \eqref{eq:Phi_convex}, denoting
				\(
					\bm v^k
				{}\coloneqq{}
					\nabla \h*(\bm x^k)-\nabla \h*(\bm u^k)
				{}\in{}
					\hat\partial\Phi(\bm u^k)
				\)
				we have
				\begin{align*}
					\FBE(\bm x^k)
					{}-{}
					\min \Phi
				{}={} &
					\Phi(\bm u^k)
					{}+{}
					\D*(\bm u^k,\bm x^k)
					{}-{}
					\Phi(\bm x^\star)
				\\
					\dueto{\cref{thm:3P}}
				{}={} &
					\Phi(\bm u^k)
					{}+{}
					\D*(\bm x^\star,\bm x^k)
					{}-{}
					\D*(\bm x^\star,\bm u^k)
					{}+{}
					\innprod*{\bm v^k}{\bm x^\star-\bm u^k}
					{}-{}
					\Phi(\bm x^\star)
				\\
					\dueto{\cref{thm:equivP:phi}}
				{}={} &
					\textstyle
					\varphi(u^k)
					{}+{}
					\D*(\bm x^\star,\bm x^k)
					{}-{}
					\D*(\bm x^\star,\bm u^k)
					{}+{}
					\innprod*{\sum_{i=1}^Nv_i^k}{x^\star-u^k}
					{}-{}
					\varphi(x^\star)
				\\
					\dueto{\cref{thm:equivP:partial}}
				{}\leq{} &
					\textstyle
					\D*(\bm x^\star,\bm x^k)
					{}-{}
					\D*(\bm x^\star,\bm u^k)
					{}-{}
					\tfrac{\mu_{\varphi}}{2}\|u^k-x^\star\|^2
				\\
					\dueto{\cref{thm:3P}}
				{}={} &
					\D*(\bm u^k,\bm x^k)
					{}+{}
					\innprod*{\nabla\h*(\bm u^k)-\nabla\h*(\bm x^k)}{\bm x^\star-\bm u^k}
					{}-{}
					\tfrac{\mu_\varphi}{2}\|u^k-x^\star\|^2
				\\
				\numberthis\label{eq:RS:strconvex:1}
				{}={} &
				\textstyle
					\D*(\bm u^k,\bm x^k)
					{}+{}
					\sum_{i=1}^N\innprod*{\nabla\h*_i(u^k)-\nabla\h*_i(x_i^k)}{x^\star-u^k}
					{}-{}
					\tfrac{\mu_\varphi}{2}\|u^k-x^\star\|^2,
				\end{align*}
				where the inequality follows from strong convexity of $\varphi$. For any \(\varepsilon_i>0\), \(i\in[N]\), one has
				\begin{align*}
					\innprod*{\nabla\h*_i(u^k)-\nabla\h*_i(x_i^k)}{x^\star-u^k}
				{}\leq{} &
					\tfrac{\varepsilon_i}{2}
					\|x^\star-u^k\|^2
					{}+{}
					\tfrac{1}{2\varepsilon_i}
					\|\nabla\h*_i(u^k)-\nabla\h*_i(x_i^k)\|^2.
				\end{align*}
				Plugged in \eqref{eq:RS:strconvex:1} with $\varepsilon_i>0$ such that $\sum_{i=1}^N\varepsilon_i= \mu_\varphi$,
				so as to cancel the square norm therein,%
				\begin{align*}
					\FBE(\bm x^k)
					{}-{}
					\min \Phi
				{}\leq{} &
					\D*(\bm u^k,\bm x^k)
					{}+{}
					\textstyle\sum_{i=1}^N{
						\tfrac{1}{2\varepsilon_i}
						\|\nabla\h*_i(u^k)-\nabla\h*_i(x_i^k)\|^2
					}
				\\
					\dueto{\cref{thm:hC11}}
				{}\leq{} &
					\D*(\bm u^k,\bm x^k)
					{}+{}
					\textstyle\sum_{i=1}^N{
						\tfrac{\ell_{\h*_i,\U}}{\varepsilon_i}
						\D*_i(u^k,x_i^k)
					}
				{}={}
					\textstyle\sum_{i=1}^N{
						\left(
							1+\tfrac{\ell_{\h*_i,\U}}{\varepsilon_i}
						\right)
						\D*_i(u^k,x_i^k)
					},
				\end{align*}
				where \(\ell_{\h*_i,\U}\) is a Lipschitz constant for \(\nabla\h*_i\) on \(\U\) as in \cref{thm:Legendre:smooth}.
				By choosing
				\(
					\varepsilon_i
				{}={}
					\nicefrac{
						\ell_{\h*_i,\U}
					}{
						\kappa
					}
				\)
				with
				\(
					\kappa
				{}\coloneqq{}
					\tfrac{
						\sum_j\ell_{\h*_j,\U}
					}{
						\mu_\varphi
					}
				\)
				(which satisfies $\sum_{i=1}^N\varepsilon_i=\mu_\varphi$), we obtain
				\[\textstyle
					\FBE(\bm x^k)
					{}-{}
					\min \Phi
				{}\leq{}
					(1+\kappa)
					\sum_{i=1}^N{
						\D*_i(u^k,x_i^k)
					}
				{}={}
					(1+\kappa)
					\D*(\bm u^k,\bm x^k).
				\]
				Combining this with \eqref{eq:RS:descent} (recall the equivalences in \cref{thm:equivBC}) yields
				\[
					\E{\FBE(\bm x^{k+1})-\min\Phi}
				{}\leq{}
					\left(
						1
						{}-{}
						\tfrac{p_{\rm min}}{1+\kappa}
					\right)
					\bigl(\FBE(\bm x^k)-\min\Phi\bigr)
				{}={}
					(1-c_{\U})
					\bigl(\FBE(\bm x^k)-\min\Phi\bigr),
				\]
				where \(c_{\U}\) as in the statement is obtained by using the estimates of \cref{thm:Legendre:smooth} for the moduli \(\ell_{\h*_i,\U}\) appearing in the constant \(\kappa\)
				(since \(\h*_i=\nicefrac{\h_i}{\gamma_i}-\nicefrac{f_i}N\) and \(\ell_{f_i,\U}\coloneqq\ell_{\h_i,\U}L_{f_i}\) is a Lipschitz modulus for \(\nabla f_i\) on \(\U\) by \cref{thm:fC11}, one has $\sigma_{f_i,\U}\geq -\ell_{f_i,\U}$ and
				\(
					\ell_{\h*_i,\U}
				{}\leq{}
					\nicefrac{\ell_{\h_i,\U}}{\gamma_i}
					{}-{}
					\nicefrac{\sigma_{f_i,\U}}N
				\)).
				This concludes the proof of \eqref{eq:RS:strconvex:Qlinear}.
				In turn, \eqref{eq:RS:strconvex:Rlinear} follows by taking unconditional expectation and using the fact that
				\(
					\varphi(z^k)
				{}={}
					\Phi(\bm u^k)
				{}\leq{}
					\FBE(\bm x^k)
				\),
				owing to \cref{thm:equivBC:phi,thm:leq}.
			\end{appendixproof}
			
			\begin{appendixproof}{thm:ECS:subseq}%
				~\par\noindent
				\bl{We use} the simpler setting of \cref{alg:BCPP} owing to the equivalence between the algorithms shown in \cref{thm:equivBC}.
				\bl*{%
					Note that if \(\varphi\) is level bounded, then by \cref{thm:bounded} a bounded convex set \(\bm\U\) exists that contains \(\seq{\bm x^k}\) and \(\seq{\bm u^k}\).
					Local Lipschitz differentiability and local strong convexity thus imply through \cref{thm:Legendre:strconvex,thm:Legendre:smooth} that \(\h*\) is \(\mu_{\h*,\bm\U}\)-strongly convex and \(\ell_{\h*,\bm\U}\)-Lipschitz differentiable on \(\bm \U\) for some constants \(\ell_{\h*,\bm\U}\geq\mu_{\h*,\bm\U}>0\).
					If, instead, those properties hold globally, then the same claims can hold with \(\bm\U=\R^{nN}\).
					Therefore, since \(g\) is assumed to be convex, either one among \cref{ass:ECS:subseq:global,ass:ECS:subseq:local} is enough to invoke \cref{thm:TLip}, implying that \(\prox_\Phi^{\h*}\) is \(\lambda\)-Lipschitz continuous on \(\bm\U\) for some \(\lambda>0\).
				}%
				
				Since all indices are updated at least once every \(T\) iterations, one has that
				\[
					\ki
				{}\coloneqq{}
					\min\set{t\in[T]}[
						\text{\(i\) is sampled at iteration \(\nu T+t-1\)}
					]
				\]
				is well defined for each index \(i\in[N]\) and \(\nu\in\N\).
				In other words, since \(i\) is sampled at iteration \(\nu T+\ki-1\) and not in any one between \(\nu T\) and \(\nu T+\ki-2\), it holds that
				\begin{equation}\label{eq:ECS:xi}
					x_i^{\nu T}
				{}={}
					x_i^{\nu T+1}=\dots=x_i^{\nu T+\ki-1}
				\quad\text{and}\quad
					x_i^{\nu T+\ki}
				{}={}
					u^{\nu T+\ki-1}
				\quad\forall i\in[N],\nu\in\N,
				\end{equation}
				recalling $\bm u^k =(u^k,\ldots, u^k)$. We now proceed to establish a descent inequality for \(\FBE\) holding every interval of \(T\) iterations.
				First,
				\begin{align*}
					\FBE(\bm x^{T(\nu+1)})
					{}-{}
					\FBE(\bm x^{\nu T})
				{}={} &
					\textstyle
					\sum_{\tau=1}^T\left(
						\FBE(\bm x^{\nu T+\tau})
						{}-{}
						\FBE(\bm x^{\nu T+\tau-1})
					\right)
				{}\overrel*[\leq]{\ref{thm:Igeq}}
					\textstyle
					-\sum_{\tau=1}^T{
						\D*(\bm x^{\nu T+\tau},\bm x^{\nu T+\tau-1})
					}
				\\
				{}\leq{} &
					-\D*(\bm x^{\nu T+t},\bm x^{\nu T+t-1})
				{}\overrel[\leq]{\ref{thm:hstrconvex}}{}
					-\tfrac{\mu_{\h*,\bm\U}}{2}
					\|\bm x^{\nu T+t}-\bm x^{\nu T+t-1}\|^2
				\numberthis\label{eq:ECS:SD1}
				\end{align*}
				holds for all \(t\in[T]\).
				Next, for every \(i\in[N]\) it holds that
				\begin{align}
				\nonumber
					\|u^{\nu T+\ki-1}-u^{\nu T}\|
				{}={} &
					\tfrac{1}{\sqrt N}
					\|\bm u^{\nu T+\ki-1}-\bm u^{\nu T}\|
				{}\leq{}
					\tfrac{\lambda}{\sqrt N}
					\|\bm x^{\nu T+\ki-1}-\bm x^{\nu T}\|
				\\
				\label{eq:ECS:uki}
				{}\leq{} &
					\tfrac{\lambda}{\sqrt N}
					\smash{\textstyle\sum_{\tau=1}^{\ki-1}\|\bm x^{\nu T+\tau}-\bm x^{\nu T+\tau-1}\|}
				{}\overrel[\leq]{\eqref{eq:ECS:SD1}}{}
					\tfrac{\lambda T}{\sqrt N}
					\sqrt{
						\tfrac{2}{\mu_{\h*,\bm\U}}
						\bigl(
							\FBE(\bm x^{\nu T})
							{}-{}
							\FBE(\bm x^{T(\nu+1)})
						\bigr)
					},
				\end{align}
				where the first inequality uses the \(\lambda\)-Lipschitz continuity of \(\prox_\Phi^{\h*}\), the second one the triangular inequality, and the last one the fact that \(\ki\leq T\).
				For all \(i\in[N]\), it follows from \cref{thm:hC11} and the triangular inequality that
				\begin{align*}
					\|x_i^{\nu T}-u^{\nu T}\|
				{}\leq{} &
					\|x_i^{\nu T}-u^{\nu T+\ki-1}\|
					{}+{}
					\|u^{\nu T+\ki-1}-u^{\nu T}\|
				\\
					\dueto{\eqref{eq:ECS:xi}}
				{}={} &
					\|x_i^{\nu T+\ki-1}-x_i^{\nu T+\ki}\|
					{}+{}
					\|u^{\nu T+\ki-1}-u^{\nu T}\|
				\\
				\numberthis\label{eq:sqrtbound}
					\dueto{\eqref{eq:ECS:SD1}~\eqref{eq:ECS:uki}}
				{}\leq{} &
					\sqrt{\tfrac{2}{\mu_{\h*,\bm\U}}}
					\Bigl(
						1+\tfrac{\lambda T}{\sqrt N}
					\Bigr)
					\sqrt{
						\FBE(\bm x^{\nu T})
						{}-{}
						\FBE(\bm x^{T(\nu+1)})
					}.
				\end{align*}
				By squaring and summing over \(i\in[N]\) we obtain
				\begin{equation}\label{eq:ECS:SD}
					\FBE(\bm x^{T(\nu+1)})
					{}-{}
					\FBE(\bm x^{\nu T})
				{}\leq{}
					-\tfrac{
						\mu_{\h*,\bm\U}
					}{
						2(1\,+\,\lambda T/\!\!\sqrt N\,\,)^2
					}
					\|\bm x^{\nu T}-\bm u^{\nu T}\|^2.
				\end{equation}
				\bl*{%
					Since by \cref{eq:ECS} in any interval of length \(T\) every index is updated at least once, by suitably shifting, for every \(t\in [T]\) the same holds for the sequences \(\seq{\bm x^{\nu T + t}}[\nu \in \N]\) and \(\seq{\bm u^{\nu T + t}}[\nu \in \N]\).
					Thus,
				}%
				\begin{align}
					\FBE(\bm x^{k+T})
					{}-{}
					\FBE(\bm x^k)
				&{}\leq{}
					-\tfrac{
						\mu_{\h*,\bm\U}
					}{
						2(1\,+\,\lambda T/\!\!\sqrt N\,\,)^2
					}
					\|\bm x^k-\bm u^k\|^2\label{eq:ECS:SDt}
					\\
				\dueto{\cref{thm:hC11}}	& {}\leq{}
					-\tfrac{
						\mu_{\h*,\bm\U}
					}{
						L_{\h*,\bm\U}(1\,+\,\lambda T/\!\!\sqrt N\,\,)^2
					}
					\D*(\bm u^k,\bm x^k)
					\quad\forall k\in\N.
				\nonumber
				\end{align}
				By telescoping the inequality and using the fact that the envelope is lower bounded (\cref{lem:env:prop:inf} and \cref{ass:argmin}), all the assertions follow from \cref{thm:subseq}.
			\end{appendixproof}
			
			\begin{appendixproof}{thm:equivLM}%
				Note that, in \cref{alg:Finito}, $\sum_{i=1}^N s_i^0 = \s^0$.
				By induction, suppose that $\sum_{i=1}^N s_i^k = \s^k$ for some $k\geq 0$.
				Then,%
				\begin{equation}\label{eq:LM:ind_si_s}
				\textstyle
					\s^{k+1}
				{}={}
					\s^k + \sum_{i\in \mathcal I^{k+1}} (s_i^{k+1}-s_i^k)
				{}={}
					\sum_{i=1}^N s_i^k + \sum_{i\in \mathcal I^{k+1}} (s_i^{k+1}-s_i^k)
				{}={}
					\sum_{i=1}^N s_i^{k+1},
				\end{equation}
				where the first equality follows by \cref{state:Finito:s+}, the second one from the induction hypothesis, and the last one from the fact that $s_i^{k+1}=s_i^k$ for $i\notin\mathcal I^{k+1}$ as in \cref{state:Finito:s+}.
				
				In what follows, let
				\(
					\mathcal N_{\rm full}
				{}\coloneqq{}
					\set{k}[\mathcal K^k=\emptyset]
				\),
				and observe that \(k\in\mathcal N_{\rm full}\) iff the if statement at \cref{state:LM:full} is true.
				In particular, it follows from \cref{state:LM:IK+_full} that
				\begin{equation}\label{eq:LM:I+_full}
					\mathcal N_{\rm full}
				{}\subseteq{}
					\set{k}[\mathcal I_\LM^{k+1}={[N]}].
				\end{equation}
				We now proceed by induction on \(k\) to establish that \(\seq{z^k_\LM}\) and \(\seq{\s^k_\LM}\) are sequences generated by \cref{alg:Finito} with index sets being chosen as
				\begin{equation}\label{eq:LM:I=I}
					\mathcal I^{k+1}
				{}\coloneqq{}
					\mathcal I_\LM^{k+1}
				\quad
					\forall k\in\N.
				\end{equation}
				The claim is true for \(k=0\); suppose it holds up to iteration \(k\geq0\).
				We consider two cases:%
				\begin{proofitemize}
				\item{\it Case 1: \(k\in\mathcal N_{\rm full}\)}.
					We have
					\begin{equation}\label{eq:LM:kinI}
					\textstyle	
						\s_\LM^{k+1}
					{}={}
						\sum_{i=1}^N\nabla\h*_i(z^{k}_\LM)
					{}={}
						\sum_{i=1}^N\nabla\h*_i(z^{k}) 
					{}={}
						\sum_{i=1}^N s_i^{k+1}
					{}={}
						\s^{k+1},
					\end{equation}
					where the first equality follows by \cref{state:LM:ts+_full}, the second one by induction, the third one by the fact that \(\mathcal I^{k+1}=\mathcal I_\LM^{k+1}=[N]\) (cf. \eqref{eq:LM:I+_full} and \eqref{eq:LM:I=I}), and the last one from \eqref{eq:LM:ind_si_s}.
					It follows that the minimization problems defining \(z^{k+1}\) and \(z_\LM^{k+1}\) (at \cref{state:Finito:z} and \cref{state:LM:z}, respectively) coincide, thus ensuring that \(z^{k+1}=z_\LM^{k+1}\) is a feasible update for \cref{alg:Finito}.
				\item{\it Case 2: \(k\notin\mathcal N_{\rm full}\)}.
					Let \(t(k)\coloneqq\max\set{t\leq k}[t\in\mathcal N_{\rm full}]\) be the last iteration before \(k\) at which the condition at \cref{state:LM:full} holds, so that, according to \cref{state:LM:tz_full,state:LM:tz}, \(\tilde z^k_\LM=z^{t(k)}_\LM\).
					We have%
					\begin{alignat*}{3}
						\s_\LM^{k+1}
					{}={} &
						\s^k_\LM
						{}+{}
						\sum_{\mathclap{i\in\mathcal I_\LM^{k+1}}}\left[
							\nabla \h*_i(z^k_\LM)
							{}-{}
							\nabla \h*_i(\tilde z^k_\LM)
						\right]
					&\qquad
						\dueto{\cref{state:LM:ts+}}
					\\
					{}={} &
						\s_{\phantom\LM}^k
						{}+{}
						\sum_{\mathclap{i\in\mathcal I^{k+1}}}\left[
							\nabla\h*_i(z^k)
							{}-{}
							\nabla\h*_i(z^{t(k)})
						\right]
					&\qquad
						\dueto{induction and \eqref{eq:LM:I=I}}
					\\
					{}={} &
						\s_{\phantom\LM}^k
						{}+{}
						\sum_{\mathclap{i\in\mathcal I^{k+1}}}\left[
							s_i^{k+1}
							{}-{}
							s_i^{t(k)+1}
						\right],
					&\qquad
						\dueto{\cref{state:Finito:s+}, \eqref{eq:LM:I+_full}, and \eqref{eq:LM:I=I}}
					\numberthis\label{eq:LM:s/2}
					\end{alignat*}
					where \eqref{eq:LM:I+_full} was used to infer that \(s_i^{t(k)+1}=\nabla\h*_i(z^{t(k)})\) for all \(i\in[N]\) (hence for all \(i\in\mathcal I^{k+1}\)).
					To conclude, note that the selection rule for \(\mathcal I_\LM^{k+1}\) (cf. \cref{state:LM:I+,state:LM:K+}) ensures through \eqref{eq:LM:I=I} that the index sets \(\mathcal I^{t(k)+1},\mathcal I^{t(k)+2},\dots,\mathcal I^{k+1}\) are all pairwise disjoint, hence that \(s_i^{t(k)}=s_i^{t(k)+1}=\dots=s_i^k\) for all \(i\in\mathcal I^{k+1}\), as is apparent from \cref{state:Finito:s+}.
					We may thus replace \(s_i^{t(k)+1}\) with \(s_i^k\) in \eqref{eq:LM:s/2} to obtain the \(\s\)-update of \cref{state:Finito:s+}, and conclude that \(\s_\LM^{k+1}=\s^{k+1}\).
					As discussed in the last part of {\it Case 1}, this in turn shows that \(z^{k+1}=z_\LM^{k+1}\) is a feasible update for \cref{alg:Finito}.
					\qedhere
				\end{proofitemize}
			\end{appendixproof}
	}%
	\end{appendix}%


	\bibliographystyle{plain}%
	\bibliography{Bibliography.bib}%

\end{document}